\documentclass[12pt,oneside,reqno]{amsart}
\usepackage{mathrsfs}
\usepackage{graphics}
\usepackage{amssymb}
\usepackage{enumerate}
\newcommand{\dif}{\mathrm{d}}

\newcommand{\be}{\begin{eqnarray}}
\newcommand{\ee}{\end{eqnarray}}
\newcommand{\ce}{\begin{eqnarray*}}
\newcommand{\de}{\end{eqnarray*}}
\newtheorem{theorem}{Theorem}[section]
\newtheorem{lemma}[theorem]{Lemma}
\newtheorem{remark}[theorem]{Remark}
\newtheorem{definition}[theorem]{Definition}
\newtheorem{proposition}[theorem]{Proposition}
\newtheorem{Example}[theorem]{Example}
\newtheorem{corollary}[theorem]{Corollary}
\def\e{\varepsilon}

\def\[{{\Big[}}
\def\]{{\Big]}}
\def\<{{\langle}}
\def\>{{\rangle}}
\def\({{\Big(}}
\def\){{\Big)}}

\def\no{\nonumber}
\def\bt{\begin{theorem}}
\def\et{\end{theorem}}
\def\bl{\begin{lemma}}
\def\el{\end{lemma}}
\def\br{\begin{remark}}
\def\er{\end{remark}}
\def\bx{\begin{Example}}
\def\ex{\end{Example}}
\def\bd{\begin{definition}}
\def\ed{\end{definition}}
\def\bp{\begin{proposition}}
\def\ep{\end{proposition}}
\def\bc{\begin{corollary}}
\def\ec{\end{corollary}}

\def\cB{{\mathcal B}}
\def\cC{{\mathcal C}}

\def\cL{{\mathcal L}}
\def\cM{{\mathcal M}}

\def\cP{{\mathcal P}}

\def\mE{{\mathbb E}}

\def\mH{{\mathbb H}}

\def\mP{{\mathbb P}}

\def\mR{{\mathbb R}}

\def\mU{{\mathbb U}}

\def\geq{\geqslant}
\def\leq{\leqslant}

\begin{document}

\allowdisplaybreaks

\title{Nonlinear filtering of stochastic differential equations with correlated L\'evy noises*}

\author{Huijie Qiao}

\dedicatory{School of Mathematics,
Southeast University\\
Nanjing, Jiangsu 211189,  China\\
Department of Mathematics, University of Illinois at
Urbana-Champaign\\
Urbana, IL 61801, USA\\
hjqiaogean@seu.edu.cn}

\thanks{{\it AMS Subject Classification(2010):} 60G57; 60G35; 60H15}

\thanks{{\it Keywords:} Nonlinear filtering problems; correlated L\'evy noises; the Kushner-Stratonovich and Zakai equations; the pathwise uniqueness and uniqueness in joint law}

\thanks{*This work was partly supported by NSF of China (No. 11001051, 11371352) and China Scholarship Council under Grant No. 201906095034.}

\subjclass{}

\date{}

\begin{abstract}
The work concerns nonlinear filtering problems of stochastic differential equations with correlated L\'evy noises. First, we establish the Kushner-Stratonovich and Zakai equations through martingale representation theorems and the Kallianpur-Striebel formula. Second, we show the pathwise uniqueness and uniqueness in joint law of weak solutions for the Zakai equation. Finally, we investigate the uniqueness in joint law of weak solutions to the Kushner-Stratonovich equation.
\end{abstract}

\maketitle \rm

\section{Introduction}

A nonlinear filtering problem means that, given a partial observable system $(X_{\cdot}, Y_{\cdot})$ defined on a complete filtered probability space $(\Omega, \mathscr{F}, \{\mathscr{F}_t\}_{t\in[0,T]},\mP)$ for $T>0$, we estimate $X_{\cdot}$ by $Y_{\cdot}$. As usual, the process $X_{\cdot}$ is difficult to observe, while the process $Y_{\cdot}$ is easy to observe and contains information about $X_{\cdot}$. The nonlinear filtering of $X_t$ with respect to $\{Y_t, 0\leq t\leq T\}$ is the `filter' $\mE[F(X_t)|\mathscr{F}_t^Y]$, where $\mathscr{F}_t^Y$ is the $\sigma$-algebra generated by $\{Y_s, 0\leq s\leq t\}$ and $F$ is any Borel measurable function with $\mE|F(X_t)|<\infty$ for $t\in[0,T]$. 

Nowadays, nonlinear filtering problems have been widely applied to various fields, such as physics, biology, the control theory and the weather forecast. Moreover, more and more researchers are paying attention to nonlinear filtering problems. At present, there have been many results about nonlinear filtering problems of stochastic systems with Gaussion noises. However, a lot of phenomena can only be simulated by stochastic systems with L\'evy noises. Therefore, nonlinear filtering problems of stochastic systems with L\'evy noises are highlighted. If the driving noises of $X_{\cdot}$ are independent of that for $Y_{\cdot}$, the type of nonlinear filtering problems has been studied in \cite{mmbp, mbfp, spss, qd}. In the paper, we solve nonlinear filtering problems of $X_t$ with respect to $\{Y_t, 0\leq t\leq T\}$ where the driving noises of $X_{\cdot}$ are correlated with that for $Y_{\cdot}$. 

Here, we explain our nonlinear filtering problems in details. Let $B, W$ be $d$-dimensional and $m$-dimensional Brownian motions defined on $(\Omega, \mathscr{F}, \{\mathscr{F}_t\}_{t\in[0,T]},\mP)$, respectively. Besides, let ($\mU,\mathscr{U}$) be a finite dimensional, measurable normed space  with the norm $\|\cdot\|_{\mU}$. And let $\nu_1$ be a $\sigma$-finite measure defined on it. Based on \cite[Theorem 9.1, P. 44]{iw}, we know that there eixsts a stationary Poisson point process $p_{\cdot}$ of the class (quasi left-continuous) with values in $\mU$ and the characteristic measure $\nu_1$ (See \cite[P. 43]{iw} for the definition of stationary Poisson point processes and \cite[P. 59]{iw} for the definition of the class (quasi left-continuous)). Let $N_{p}((0,t],\dif u)$ be the counting measure of $(p_{t})$ such that $\mE N_{p}((0,t],A)=t\nu_1(A)$ for $A\in\mathscr{U}$. Denote
$$
\tilde{N}_{p}((0,t],\dif u):=N_{p}((0,t],\dif u)-t\nu_1(\dif u),
$$
the compensated martingale measure of $p_{t}$. Fix $\mU_1\in\mathscr{U}$ with $\nu_1(\mU\setminus\mU_1)<\infty$ and $\int_{\mU_1}\|u\|_{\mU}^2\,\nu_1(\dif u)<\infty$. Let $X_{\cdot}$ be the solution process of the following stochastic differential equation (SDE in short) on $\mR^n$:
\ce
\dif X_t=b_1(t,X_t)\dif t+\sigma_0(t,X_t)\dif B_t+\sigma_1(t,X_t)\dif W_t+\int_{\mU_1}f_1(t,X_{t-},u)\tilde{N_p}(\dif t, \dif u), \quad 0\leq t\leq T,
\de
where the mappings $b_1:[0,T]\times\mR^n\mapsto\mR^n$, $\sigma_0:[0,T]\times\mR^n\mapsto\mR^{n\times d}$, $\sigma_1:[0,T]\times\mR^n\mapsto\mR^{n\times m}$ and $f_1:[0,T]\times\mR^n\times\mU_1\mapsto\mR^n$ are all Borel measurable. And then $Y_{\cdot}$ is the solution process of the following SDE on $\mR^m$:
\ce
\dif Y_t=b_2(t,X_t,Y_t)\dif t+\sigma_2(t,Y_t)\dif W_t+\int_{\mU_2}f_2(t, Y_{t-}, u)\tilde{N}_{\lambda}(\dif t, \dif u), \quad 0\leq t\leq T,
\de
where the mappings $b_2:[0,T]\times\mR^n\times\mR^m\mapsto\mR^m$, $\sigma_2:[0,T]\times\mR^m\mapsto\mR^{m\times m}$, $f_2:[0,T]\times\mR^m\times\mU_2\mapsto\mR^m$
are all Borel measurable, and $N_{\lambda}(\dif t,\dif u)$ is an integer-valued random measure and its predictable compensator is given by $\lambda(t,X_{t-},u)\dif t\nu_2(\dif u)$. That is, $\tilde{N}_\lambda(\dif t,\dif u):=N_\lambda(\dif t,\dif u)-\lambda(t,X_{t-},u)\dif t\nu_2(\dif u)$ is its compensated martingale measure. Here $\nu_2$ is another $\sigma$-finite measure defined on $\mU$ with $\nu_2(\mU\setminus\mU_2)<\infty$ and $\int_{\mU_2}\|u\|_{\mU}^2\,\nu_2(\dif u)<\infty$ for $\mU_2\in\mathscr{U}$. And $\lambda(t,x,u):[0,T]\times\mR^n\times\mU\mapsto(0,1)$ is Borel measurable. We will study the nonlinear filtering problem of $X_t$ with respect to $\{Y_t, 0\leq t\leq T\}$.

Our motivation is two-folded. One fold lies in these models themselves. These models are usually called as feedback models. That is, the observation $(Y_t)$ is fed back to the dynamics of the signal $(X_t)$. And feedback models have appeared in many applications (especially in aerospace problems). Note that our models are different from ones in \cite{cckc1, cckc2}, where $\sigma_1=0$ and $B_t$ and $W_t$ are correlated each other. And our models are more natural and applicable than that in \cite{cckc1, cckc2}. The other fold is that the filtering of these models can be used to solve nonlinear filtering problems of stochastic multi-scale systems with correlated L\'evy noises (\cite{q0}). And then we not only deduce the Kushner-Stratonovich and Zakai equations but also investigate the pathwise uniqueness and uniqueness in joint law of weak solutions for the two equations.

It is worthwhile to mentioning our methods. First of all, since the driving processes of $X_{\cdot}$ are not independent of $Y_{\cdot}$, the method of measure transformations does {\it not} work. Therefore, we make use of martingale representation theorems and the Kallianpur-Striebel formula to obtain the Kushner-Stratonovich and Zakai equations. The difficulty of deduction lies in looking for suitable martingales. Moreover, our method can be applied to solve nonlinear filtering problems of stochastic differential equations with correlated sensor L\'evy noises (See Section \ref{con} for details). Second, about the uniqueness of solutions for the two equations, there are two methods--the method of filtered martingale problems (\cite{cckc1, qd}) and the method of operator equations (\cite{q2}). Specially, in \cite{qd}, Qiao and Duan required that $\lambda$ in the driving processes of the observation process $Y_{\cdot}$  is independent of $x$. And in \cite{q2}   the author assumed that the driving processes of $X_{\cdot}$ has no jump term. Here we prove two types of uniquenesses for weak solutions to the two equations by a family of operators without any assumption on driving processes. Therefore, our result covers that in \cite{q2, qd}. Besides, note that in \cite{mpx} under the assumption that the driving processes of $X_{\cdot}$ are independent of $Y_{\cdot}$, Maroulas et al. showed the uniquenesses of strong solutions for the two equations by the similar method. Here, the driving processes of $X_{\cdot}$ correlate with that of $Y_{\cdot}$. Thus, our result generalizes that in \cite{mpx}. However, they did not give out the clear definitions of solutions and uniquenesses about solutions for the two equations. Since strong solutions and the pathwise uniqueness of strong solutions for the Kushner-Stratonovich equation can not be defined, the results related with them are wrong.(See Remark \ref{notun}) 

The paper is arranged as follows. In Section \ref{nonfilter}, we deduce the Kushner-Stratonovich and Zakai equations by martingale problems and the Kallianpur-Striebel formula. The pathwise uniqueness and uniqueness in joint law of weak solutions for the Zakai equation are placed in Section \ref{unzaks}. In Section \ref{unks}, we investigate the uniqueness in joint law of weak solutions for the Kushner-Stratonovich equation. And then in Section \ref{con}, we summarize all the results and point out other models where our method can be applied. Finally, we prove Lemma \ref{brmoposs}, Lemma \ref{trex} and Lemma \ref{prozak} in the appendix.

The following convention will be used throughout the paper: $C$ with or without indices will denote
different positive constants whose values may change from one place to another.

\section{Nonlinear filtering problems}\label{nonfilter}

In this section, we introduce the nonlinear filtering problem for a non-Gaussian signal-observation system with correlated noises,  and derive the Kushner-Stratonovich and Zakai equations.

\subsection{The framework}\label{fram}

In the subsection, we introduce signal-observation systems.

Consider the following signal-observation system $(X_t, Y_t)$ on $\mR^n\times\mR^m$:
\be\left\{\begin{array}{l}
\dif X_t=b_1(t,X_t)\dif t+\sigma_0(t,X_t)\dif B_t+\sigma_1(t,X_t)\dif W_t+\int_{\mU_1}f_1(t,X_{t-},u)\tilde{N_p}(\dif t, \dif u),\\
\dif Y_t=b_2(t,X_t,Y_t)\dif t+\sigma_2(t,Y_t)\dif W_t+\int_{\mU_2}f_2(t, Y_{t-}, u)\tilde{N}_{\lambda}(\dif t, \dif u), \quad 0\leq t\leq T.\\
 \end{array}
\right. 
\label{Eq1} 
\ee
 The initial value $X_0$ is assumed to be a square integrable random variable independent of $Y_0, B_{\cdot}, W_{\cdot}, N_p, N_{\lambda}$. Moreover, $B_{\cdot}, W_{\cdot}, N_p, N_{\lambda}$ are mutually independent.
We make the following hypotheses:
\begin{enumerate}[($\mathbf{H}^1_{b_1, \sigma_0, \sigma_1, f_1}$)] 
\item For $t\in[0,T]$ and $x_1, x_2\in\mR^n$,
\ce
&|b_1(t,x_1)-b_1(t,x_2)|\leq L_1(t)|x_1-x_2|\kappa_1(|x_1-x_2|),\\
&\|\sigma_0(t,x_1)-\sigma_0(t,x_2)\|^2\leq L_1(t)|x_1-x_2|^{2}\kappa_2(|x_1-x_2|),\\
&\|\sigma_1(t,x_1)-\sigma_1(t,x_2)\|^2\leq L_1(t)|x_1-x_2|^{2}\kappa_3(|x_1-x_2|),\\
&\int_{\mU_1}|f_1(t,x_1,u)-f_1(t,x_2,u)|^{p'}\,\nu_1(\dif u)\leq L_1(t)|x_1-x_2|^{p'}\kappa_4(|x_1-x_2|),
\de
hold for $p'=2$ and $4$, where $|\cdot|$ and $\|\cdot\|$ denote the Hilbert-Schmidt norms of a vector and a matrix, respectively. Here $L_1(t)>0$ is an increasing function and $\kappa_i$ is a positive continuous
function, bounded on $[1,\infty)$ and satisfies
\ce
\lim\limits_{x\downarrow0}\frac{\kappa_i(x)}{\log x^{-1}}<\infty, \quad i=1, 2, 3, 4.
\de
\end{enumerate}
\begin{enumerate}[($\mathbf{H}^2_{b_1, \sigma_0,\sigma_1, f_1}$)]
\item For $t\in[0,T]$ and $x\in\mR^n$,
$$
|b_1(t,x)|^2+\|\sigma_0(t,x)\|^2+\|\sigma_1(t,x)\|^2+\int_{\mU_1}|f_1(t,x,u)|^2\nu_1(\dif u)\leq K_1(t)(1+|x|)^2,
$$
where $K_1(t)>0$ is an increasing function.
\end{enumerate}
\begin{enumerate}[($\mathbf{H}^1_{\sigma_2, f_2}$)] 
\item For $t\in[0,T]$ and $y_1, y_2\in\mR^m$,
\ce
&\|\sigma_2(t,y_1)-\sigma_2(t,y_2)\|^2\leq L_2(t)|y_1-y_2|^2,\\
&\int_{\mU_2}|f_2(t,y_1,u)-f_2(t,y_2,u)|^2\,\nu_2(\dif u)\leq L_2(t)|y_1-y_2|^2,
\de
where $L_2(t)>0$ is an increasing function.
\end{enumerate}
\begin{enumerate}[($\mathbf{H}^2_{b_2, \sigma_2, f_2}$)] 
\item For $t\in[0,T], y\in\mR^m$, $\sigma_2(t,y)$ is invertible,
$$
|b_2(t,x,y)|\vee\|\sigma_2(t,0)\|\vee\|\sigma^{-1}_2(t,y)\|\leq K_2,~{for}~{all}~t\in[0,T], x\in\mR^n, y\in\mR^m,
$$
where $K_2>0$ is a constant,
and
$$
\sup\limits_{s\in[0,T]}\int_{\mU_2}|f_2(s,0,u)|^2\nu_2(\dif u)<\infty.
$$
\end{enumerate}

\br
The assumption ($\mathbf{H}^2_{b_2, \sigma_2, f_2}$) assures the well-posedness of strong solutions for the second equation in the system (\ref{Eq1}) and the usage of the following Girsanov theorem.
\er

By Theorem 1.2 in \cite{q1}, the system (\ref{Eq1}) has a pathwise unique strong solution denoted as $(X_t,Y_t)$. 

\subsection{Characterization of $\mathscr{F}_t^Y$}\label{char}

In the subsection, we describe $\mathscr{F}_t^Y$ by the Girsanov theorem.

First of all, we assume:
\begin{enumerate}[($\mathbf{H}_{\lambda}$)] 
\item There exists a function $L(u): \mU_2\rightarrow \mR^+$ satisfying $0<\iota<L(u)<\lambda(t,x,u)<1$ and 
$$
\int_{\mU_2}\frac{(1-L(u))^2}{L(u)}\nu_2(\dif u)<\infty,
$$
where $0<\iota<1$ is a constant, and
\be
\int_0^T\int_{\mU_2}\mE\lambda(s,X_s,u)\nu_2(\dif u)\dif s<\infty.
\label{intecond}
\ee 
\end{enumerate}

Set
\ce
h(t,X_t,Y_t):=\sigma_2^{-1}(t, Y_t)b_2(t,X_t,Y_t),
\de
\ce
\Lambda^{-1}_t:&=&\exp\bigg\{-\int_0^th^i(s,X_s,Y_s)\dif W^i_s-\frac{1}{2}\int_0^t
\left|h(s,X_s,Y_s)\right|^2\dif s\\
&&\quad\qquad -\int_0^t\int_{\mU_2}\log\lambda(s,X_{s-},u)\tilde{N}_{\lambda}(\dif s, \dif u)\\
&&\quad\qquad -\int_0^t\int_{\mU_2}\(1-\lambda(s,X_s,u)+\lambda(s,X_s,u)\log\lambda(s,X_s,u)\)\nu_2(\dif u)\dif s\bigg\}.
\de
Here and hereafter, we use the convention that repeated indices imply summation. By ($\mathbf{H}^2_{b_2, \sigma_2, f_2}$) and ($\mathbf{H}_{\lambda}$), we know that 
\ce
\mE\left(\int_0^T
\left|h(s,X_s,Y_s)\right|^2\dif s\right)<\infty,
\de
and
\ce
&&\mE\left(\int_0^T\int_{\mU_2}|\log\lambda(s,X_{s-},u)|^2\lambda(s,X_{s-},u)\nu_2(\dif u)\dif s\right)\leq\int_0^T\int_{\mU_2}|\log L(u)|^2\nu_2(\dif u)\dif s\\
&\leq& \int_0^T\int_{\mU_2}\frac{(1-L(u))^2}{L^2(u)}\nu_2(\dif u)\dif s
\leq\int_0^T\int_{\mU_2}\frac{(1-L(u))^2}{L(u)}\frac{1}{\iota}\nu_2(\dif u)\dif s<\infty.
\de
Thus, the definition of $\Lambda^{-1}_t$ is reasonable. Again set
$$
M_t:=-\int_0^t h^i(s,X_s,Y_s)\dif W^i_s
+\int_0^t\int_{\mU_2}\frac{1-\lambda(s,X_{s-},u)}{\lambda(s,X_{s-},u)}\tilde{N}_{\lambda}(\dif s, \dif u),
$$
and then by the similar deduction to \cite{qd}, we know that $\Lambda^{-1}_t$, the Dol\'eans-Dade
exponential of $M$, is an exponential martingale. Define a measure $\tilde{\mP}$ via
$$
\frac{\dif \tilde{\mP}}{\dif \mP}=\Lambda^{-1}_T.
$$
By the Girsanov theorem for Brownian motions and random measures(e.g.Theorem 3.17 in \cite{jjas}), one can obtain that under the measure $\tilde{\mP}$, 
\be\label{tilw}
\tilde{W}_t:=W_t+\int_0^t h(s,X_s,Y_s)\dif s
\ee
is an $(\mathscr{F}_t)$-adapted Brownian motion,
\be\label{tiln}
\tilde{N}(\dif t, \dif u):=N_\lambda(\dif t, \dif u)-\dif t\nu_2(\dif u),
\ee
is an $(\mathscr{F}_t)$-adapted Poisson martingale measure, and the system (\ref{Eq1}) becomes
\be\left\{\begin{array}{l}
\dif X_t=\tilde{b}_1(t,X_t)\dif t+\sigma_0(t,X_t)\dif B_t+\sigma_1(t,X_t)\dif \tilde{W}_t+\int_{\mU_1}f_1(t,X_{t-},u)\tilde{N}_p(\dif t, \dif u),\\
\dif Y_t=\sigma_2(t,Y_t)\dif \tilde{W}_t+\int_{\mU_2}f_2(t,Y_{t-},u)\tilde{N}(\dif t, \dif u),
\end{array}
\right. \label{Eq2}
\ee
where
\ce
\tilde{b}_1(t,X_t)=b_1(t,X_t)-\sigma_1(t,X_t)h(t,X_t,Y_t).
\de
Furthermore, the $\sigma$-algebra $\mathscr{F}_t^{Y^0}$ generated by $\{Y_s, 0\leq s\leq t\}$, can be characterized as
\ce
\mathscr{F}_t^{Y^0}=\mathscr{F}_t^{\tilde{W}}\vee\mathscr{F}_t^{\tilde{N}}\vee\mathscr{F}_0^{Y^0},
\de
where $\mathscr{F}_t^{\tilde{W}}, \mathscr{F}_t^{\tilde{N}}$ denote the $\sigma$-algebras generated by $\{\tilde{W}_s, 0\leq s\leq t\}, \{\tilde{N}((0,s], A), 0\leq s\leq t, A\in\mathscr{U}\}$, respectively (See \cite[Lemma 3.2]{qd} for details). And then $\mathscr{F}_t^Y$ denotes the usual augmentation of $\mathscr{F}_t^{Y^0}$. 

\subsection{The Kushner-Stratonovich equation}\label{kseq}

Next, set
\ce
\mP_t(F):=\mE[F(X_t)|\mathscr{F}_t^Y], \quad F\in\cB(\mR^n),
\de
where $\cB(\mR^n)$ denotes the set of all bounded measurable functions on $\mR^n$, and then $\mP_t$ is called as the nonlinear filtering of $X_t$ with respect to $\mathscr{F}_t^Y$. Moreover, the equation satisfied by $\mP_t$ is called the Kushner-Stratonovich equation. In order  to derive the Kushner-Stratonovich equation, we  need these following results.

\bl\label{brmoposs}
Under the measure $\mP$, $\bar{W}_t:=\tilde{W}_t-\int_0^t\mP_s(h(s,\cdot,Y_s))\dif s$ is an $(\mathscr{F}^Y_t)$-adapted Brownian motion and $\tilde{\bar{N}}(\dif t, \dif u)=N_\lambda(\dif t, \dif u)-\mP_{t-}\left(\lambda(t,\cdot,u)\right)\nu_2(\dif u)\dif t$ is an $(\mathscr{F}^Y_t)$-adapted martingale measure, where $\mP_{t-}\left(\lambda(t,\cdot,u)\right)\nu_2(\dif u)\dif t$ is the $(\mP, \mathscr{F}_t^Y)$-predictable projection of $N_\lambda(\dif t, \dif u)$.
\el

Although the result in the above lemma has appeared, we haven't seen its proof. Therefore, to the readers' convenience, the detailed proof is placed in the appendix.
 
 \br
 $\bar{W}$ is usually called the innovation process.
 \er
 
 The following lemma comes from \cite[Proposition 2.1]{cckc1}.
  
 \bl\label{mart1}
Suppose that $(M_t)$ is an $(\mathscr{F}_t)$-adapted local martingale. If there exists a localizing $(\mathscr{F}^Y_t)$-stopping time sequence $\{\tau_n\}$ for $(M_t)$, then $\(\mE[M_t|\mathscr{F}_t^Y]\)$ is an $(\mathscr{F}^Y_t)$-adapted local martingale.
\el

By the above lemma, it is obvious that if $(M_t)$ is an $(\mathscr{F}_t)$-adapted martingale, then $\(\mE[M_t|\mathscr{F}_t^Y]\)$ is a $(\mathscr{F}^Y_t)$-adapted martingale.

\bl\label{mart2}
Suppose that $(\phi_t)$ is a measurable process satisfying 
$$
\mE\[\int_0^T|\phi_s|\dif s\]<\infty.
$$
Then $\(\mE[\int_0^t\phi_s\dif s|\mathscr{F}_t^Y]-\int_0^t\mE[\phi_s|\mathscr{F}_s^Y]\dif s\)$ is an $(\mathscr{F}^Y_t)$-adapted martingale.
\el 

Since the proof of the above lemma is only based on the tower property of conditional expectations, we omit it. Now, it is the position to give and deduce the Kushner-Stratonovich equation.

\bt  (The Kushner-Stratonovich equation)  \label{ks}\\
For $F\in\cC_c^\infty(\mR^n)$, the Kushner-Stratonovich equation of the system (\ref{Eq1}) is given by
\be
\mP_t(F)&=&\mP_0(F)+\int_0^t\mP_s(\cL_s F)\dif s+\sum\limits_{i=1}^n\sum\limits_{l=1}^m\int_0^t\mP_s\left(\frac{\partial F(\cdot)}{\partial x_i}\sigma^{il}_1(s,\cdot)\right)\dif  \bar{W}^l_s\no\\
&&+\sum\limits_{l=1}^m\int_0^t\bigg(\mP_s\left(Fh^l(s,\cdot,Y_s)\right)-\mP_s\left(F\right)\mP_s\left(h^l(s,\cdot,Y_s)\right)\bigg)\dif \bar{W}^l_s\no\\
&&+\int_0^t\int_{\mU_2}\frac{\mP_{s-}\left(F\lambda(s,\cdot,u)\right)-\mP_{s-}\left(F\right)\mP_{s-}
\left(\lambda(s,\cdot,u)\right)}{\mP_{s-}\left(\lambda(s,\cdot,u)\right)}\tilde{\bar{N}}(\dif s, \dif u), t\in[0,T],
\label{kseq0}
\ee
where the operater $\cL_s$ is defined as
\ce
(\cL_s F)(x)&=&\frac{\partial F(x)}{\partial x_i}b^i_1(s,x)+\frac{1}{2}\frac{\partial^2F(x)}{\partial x_i\partial x_j}
\sigma_0^{ik}(s,x) \sigma_0^{kj}(s,x)+\frac{1}{2}\frac{\partial^2F(x)}{\partial x_i\partial x_j}
\sigma_1^{il}(s,x) \sigma_1^{lj}(s,x)\\
&&+\int_{\mU_1}\Big[F\big(x+f_1(s,x,u)\big)-F(x)
-\frac{\partial F(x)}{\partial x_i}f^i_1(s,x,u)\Big]\nu_1(\dif u).
\de
\et

\br
Here, we justify that three stochastic integrals in (\ref{kseq0}) are well defined. First of all, since $F\in\cC_c^\infty(\mR^n)$, ($\mathbf{H}^2_{b_1, \sigma_0,\sigma_1, f_1}$) ($\mathbf{H}^2_{b_2, \sigma_2, f_2}$) admit us to obtain that
\ce
&&\sum\limits_{l=1}^m\mE\left(\int_0^T\left|\mP_s\left(\sum\limits_{i=1}^n\frac{\partial F(\cdot)}{\partial x_i}\sigma^{il}_1(s,\cdot)\right)\right|^2\dif s\right)<\infty,\\
&&\sum\limits_{l=1}^m\mE\left(\int_0^t\left|\mP_s\left(Fh^l(s,\cdot,Y_s)\right)-\mP_s\left(F\right)\mP_s\left(h^l(s,\cdot,Y_s)\right)\right|^2\dif s\right)<\infty.
\de
And then by ($\mathbf{H}_{\lambda}$), it holds that
\ce
&&\mE\left(\int_0^T\int_{\mU_2}\left|\frac{\mP_{s-}\left(F\lambda(s,\cdot,u)\right)-\mP_{s-}\left(F\right)\mP_{s-}
\left(\lambda(s,\cdot,u)\right)}{\mP_{s-}\left(\lambda(s,\cdot,u)\right)}\right|^2\mP_{s-}\left(\lambda(s,\cdot,u)\right)\nu_2(\dif u)\dif s\right)\\
&\leq&2\mE\left(\int_0^T\int_{\mU_2}\frac{|\mP_{s-}\left(F\lambda(s,\cdot,u)\right)|^2+|\mP_{s-}\left(F\right)\mP_{s-}
\left(\lambda(s,\cdot,u)\right)|^2}{\mP_{s-}\left(\lambda(s,\cdot,u)\right)}\nu_2(\dif u)\dif s\right)\\
&\leq&4\|F\|_{\infty}\mE\left(\int_0^T\int_{\mU_2}\mP_{s-}\left(\lambda(s,\cdot,u)\right)\nu_2(\dif u)\dif s\right)\\
&=&4\|F\|_{\infty}\int_0^T\int_{\mU_2}\mE\left[\lambda(s,X_s,u)\right]\nu_2(\dif u)\dif s<\infty,
\de
where $\|F\|_{\infty}:=\sup\limits_{x\in\mR^n}|F(x)|$. Thus, three stochastic integrals in (\ref{kseq0}) are well defined.
\er

\begin{proof}
Applying the It\^o formula to $X_t$, we have
\be\label{fxito}
F(X_t)&=&F(X_0)+\int_0^t(\cL_s F)(X_s)\dif s
+\sum\limits_{i=1}^n\sum\limits_{k=1}^d\int_0^t\frac{\partial F(X_s)}{\partial x_i}\sigma^{ik}_0(s,X_s)\dif B^k_s\no\\
&&+\sum\limits_{i=1}^n\sum\limits_{l=1}^m\int_0^t\frac{\partial F(X_s)}{\partial x_i}\sigma^{il}_1(s,X_s)\dif W^l_s\no\\
&&+\int_0^t\int_{\mU_1}\Big[F\big(X_{s-}+f_1(s,X_{s-},u)\big)-F(X_{s-})\Big]\tilde{N}_p(\dif s, \dif u)\no\\
&=:&F(X_0)+\int_0^t(\cL_s F)(X_s)\dif s+M_t,
\ee
where $(M_t)$ is an $(\mathscr{F}_t)$-adapted local martingale. And then, by taking the conditional expectation with respect to $\mathscr{F}_t^Y$ on two hand sides of the above equality, one can obtain that
\ce
\mE[F(X_t)|\mathscr{F}_t^Y]&=&\mE[F(X_0)|\mathscr{F}_t^Y]+\mE\left[\int_0^t(\cL_s F)(X_s)\dif s|\mathscr{F}_t^Y\right]
+\mE[M_t|\mathscr{F}_t^Y],
\de
and furthermore
\ce
&&\mE[F(X_t)|\mathscr{F}_t^Y]-\mE[F(X_0)|\mathscr{F}_t^Y]-\int_0^t\mE\left[(\cL_s F)(X_s)|\mathscr{F}_s^Y\right]\dif s\\
&=&\mE\left[\int_0^t(\cL_s F)(X_s)\dif s|\mathscr{F}_t^Y\right]-\int_0^t\mE\left[(\cL_s F)(X_s)|\mathscr{F}_s^Y\right]\dif s
+\mE[M_t|\mathscr{F}_t^Y].
\de
Based on Lemma \ref{mart1} and \ref{mart2}, it holds that the right hand side of the above equality is an $(\mathscr{F}^Y_t)$-adapted local martingale. Thus, by Corollary III 4.27 in \cite{jjas}, there exist a $m$-dimensional $(\mathscr{F}^Y_t)$-adapted process $(E_t)$ and a $1$-dimensional $(\mathscr{F}^Y_t)$-predictable process $\(D(t,u)\)$ such that 
\ce
\mE[F(X_t)|\mathscr{F}_t^Y]-\mE[F(X_0)|\mathscr{F}_t^Y]-\int_0^t\mE\left[(\cL_s F)(X_s)|\mathscr{F}_s^Y\right]\dif s\\
=\int_0^tE_s\dif\bar{W}_s+\int_0^t\int_{\mU_2}D(s,u)\tilde{\bar{N}}(\dif s, \dif u).
\de
Note that $\mP_t(F)=\mE[F(X_t)|\mathscr{F}_t^Y]$ and $X_0$ is independent of $(\mathscr{F}^Y_t)$. And then we have that
\be
\mP_t(F)=\mP_0(F)+\int_0^t\mP_s\left(\cL_s F\right)\dif s+\int_0^tE_s\dif\bar{W}_s+\int_0^t\int_{\mU_2}D(s,u)\tilde{\bar{N}}(\dif s, \dif u).
\label{cdun}
\ee

Next let us firstly determine $(E_t)$. On one side, applying the It\^o formula to $F(X_t)\tilde{W}^j_t$, by (\ref{tilw}) (\ref{fxito}) one can get that for $j=1, 2, \cdots, m$,
\ce
F(X_t)\tilde{W}^j_t&=&\int_0^tF(X_s)\dif \tilde{W}^j_s+\int_0^t\tilde{W}^j_s\dif F(X_s)+\int_0^t\frac{\partial F(X_s)}{\partial x_i}\sigma^{ij}_1(s,X_s)\dif s\\
&=&\int_0^tF(X_s)h^j(s,X_s,Y_s)\dif s+\int_0^t\tilde{W}^j_s(\cL_s F)(X_s)\dif s+\int_0^t\frac{\partial F(X_s)}{\partial x_i}\sigma^{ij}_1(s,X_s)\dif s\\
&&+\int_0^tF(X_s)\dif W^j_s+\int_0^t\tilde{W}^j_s\dif M_s.
\de
Taking the conditional expectation with respect to $\mathscr{F}_t^Y$, by the measurability of $\tilde{W}^j_t$ with respect to $\mathscr{F}_t^Y$ we know that 
\be
\mE[F(X_t)|\mathscr{F}_t^Y]\tilde{W}^j_t&=&\int_0^t\mE\left[F(X_s)h^j(s,X_s,Y_s)|\mathscr{F}_s^Y\right]\dif s+\int_0^t\tilde{W}^j_s\mE\left[(\cL_s F)(X_s)|\mathscr{F}_s^Y\right]\dif s\no\\
&&+\int_0^t\mE\left[\frac{\partial F(X_s)}{\partial x_i}\sigma^{ij}_1(s,X_s)|\mathscr{F}_s^Y\right]\dif s+\Sigma^1_t,
\label{cont1}
\ee
where $(\Sigma^1_t)$ is an $(\mathscr{F}^Y_t)$-adapted local martingale and given by
\ce
\Sigma^1_t&:=&\left(\mE\left[\int_0^tF(X_s)h^j(s,X_s,Y_s)\dif s|\mathscr{F}_t^Y\right]-\int_0^t\mE\left[F(X_s)h^j(s,X_s,Y_s)|\mathscr{F}_s^Y\right]\dif s\right)\no\\
&&+\left(\mE\left[\int_0^t\tilde{W}^j_s(\cL_s F)(X_s)\dif s|\mathscr{F}_t^Y\right]-\int_0^t\tilde{W}^j_s\mE\left[(\cL_s F)(X_s)|\mathscr{F}_s^Y\right]\dif s\right)\no\\
&&+\left(\mE\left[\int_0^t\frac{\partial F(X_s)}{\partial x_i}\sigma^{ij}_1(s,X_s)\dif s|\mathscr{F}_t^Y\right]-\int_0^t\mE\left[\frac{\partial F(X_s)}{\partial x_i}\sigma^{ij}_1(s,X_s)|\mathscr{F}_s^Y\right]\dif s\right)\no\\
&&+\mE\left[\int_0^tF(X_s)\dif W^j_s|\mathscr{F}_t^Y\right]+\mE\left[\int_0^t\tilde{W}^j_s\dif M_s|\mathscr{F}_t^Y\right].
\de

On the other side, one can apply the It\^o formula to $\mE[F(X_t)|\mathscr{F}_t^Y]\tilde{W}^j_t$ and obtain that 
\be
 \mE[F(X_t)|\mathscr{F}_t^Y]\tilde{W}^j_t&=&\int_0^t\mE[F(X_s)|\mathscr{F}_s^Y]\dif \tilde{W}^j_s+\int_0^t\tilde{W}^j_s\dif \mE[F(X_s)|\mathscr{F}_s^Y]+\int_0^tE^j_s\dif s\no\\
 &=&\int_0^t\mE[F(X_s)|\mathscr{F}_s^Y]\mE[h^j(s,X_s,Y_s)|\mathscr{F}_s^Y]\dif s+\int_0^t\tilde{W}^j_s\mE[(\cL_s F)(X_s)|\mathscr{F}_s^Y]\dif s\no\\
&&+\int_0^tE^j_s\dif s+\Sigma^2_t,
\label{cont2}
\ee
where $(\Sigma^2_t)$ is an $(\mathscr{F}^Y_t)$-adapted local martingale and given by
\ce
\Sigma^2_t:=\int_0^t\mE[F(X_s)|\mathscr{F}_s^Y]\dif \bar{W}^j_s+\int_0^t\tilde{W}^j_sE_s\dif\bar{W}_s
+\int_0^t\int_{\mU_2}\tilde{W}^j_sD(s,u)\tilde{\bar{N}}(\dif s, \dif u).
\de

Since the left side of (\ref{cont1}) is the same to that of (\ref{cont2}), bounded variation parts of their right sides should be the same. Therefore,
\be
E^j_s&=&\mE[F(X_s)h^j(s,X_s,Y_s)|\mathscr{F}_s^Y]-\mE[F(X_s)|\mathscr{F}_s^Y]\mE[h^j(s,X_s,Y_s)|\mathscr{F}_s^Y]\no\\
&&+\mE\left[\frac{\partial F(X_s)}{\partial x_i}\sigma^{ij}_1(s,X_s)|\mathscr{F}_s^Y\right]\no\\
&=&\mP_s\(Fh^j(s,\cdot,Y_s)\)-\mP_s(F)\mP_s\(h^j(s,\cdot,Y_s)\)+\mP_s\left(\frac{\partial F(\cdot)}{\partial x_i}\sigma^{ij}_1(s,\cdot)\right), \quad a.s. \mP.
\label{cpro}
\ee

In the following we search for $\(D(t,u)\)$. Take 
$$
Z_t:=\int_0^t\int_{\mU_2}\(1-L(u)\)\tilde{N}(\dif s, \dif u), \quad 0\leq t\leq T.
$$
On one side, it follows from the It\^o formula for $F(X_t)Z_t$ that 
\ce
F(X_t)Z_t&=&\int_0^tF(X_s)\dif Z_s+\int_0^tZ_s\dif F(X_s)\\
&=&\int_0^t\int_{\mU_2}F(X_s)\(1-L(u)\)\lambda(s,X_s,u)\nu_2(\dif u)\dif s-\int_0^t\int_{\mU_2}F(X_s)\(1-L(u)\)\nu_2(\dif u)\dif s\\
&&+\int_0^tZ_s(\cL_s F)(X_s)\dif s+\int_0^t\int_{\mU_2}F(X_s)\(1-L(u)\)\tilde{N}_\lambda(\dif s, \dif u)+\int_0^tZ_s\dif M_s.
\de
Taking the conditional expectation with respect to $\mathscr{F}_t^Y$, by the measurability of $Z_t$ with respect to $\mathscr{F}_t^Y$ we get that
\be
\mE[F(X_t)|\mathscr{F}_t^Y]Z_t&=&\int_0^t\int_{\mU_2}\(1-L(u)\)\mE[F(X_s)\lambda(s,X_s,u)|\mathscr{F}_s^Y]\nu_2(\dif u)\dif s\no\\
&&-\int_0^t\int_{\mU_2}\(1-L(u)\)\mE[F(X_s)|\mathscr{F}_s^Y]\nu_2(\dif u)\dif s\no\\
&&+\int_0^tZ_s\mE[(\cL_s F)(X_s)|\mathscr{F}_s^Y]\dif s+\Sigma^3_t,
\label{disc1}
\ee
where $(\Sigma^3_t)$ is an $(\mathscr{F}^Y_t)$-adapted local martingale and given by
\ce
\Sigma^3_t&:=&\Bigg(\mE\left[\int_0^t\int_{\mU_2}F(X_s)\(1-L(u)\)\lambda(s,X_s,u)\nu_2(\dif u)\dif s|\mathscr{F}_t^Y\right]\\
&&\qquad -\int_0^t\int_{\mU_2}\(1-L(u)\)\mE[F(X_s)\lambda(s,X_s,u)|\mathscr{F}_s^Y]\nu_2(\dif u)\dif s\Bigg)\\
&&-\Bigg(\mE\left[\int_0^t\int_{\mU_2}F(X_s)\(1-L(u)\)\nu_2(\dif u)\dif s|\mathscr{F}_t^Y\right]\\
&&\qquad -\int_0^t\int_{\mU_2}\(1-L(u)\)\mE[F(X_s)|\mathscr{F}_s^Y]\nu_2(\dif u)\dif s\Bigg)\\
&&+\left(\mE\left[\int_0^tZ_s(\cL_s F)(X_s)\dif s|\mathscr{F}_t^Y\right]-\int_0^tZ_s\mE[(\cL_s F)(X_s)|\mathscr{F}_s^Y]\dif s\right)\\
&&+\mE\left[\int_0^tZ_s\dif M_s|\mathscr{F}_t^Y\right]+\mE\left[\int_0^t\int_{\mU_2}F(X_s)\(1-L(u)\)\tilde{N}_\lambda(\dif s, \dif u)|\mathscr{F}_t^Y\right].
\de

On the other side, by making use of the It\^o formula one can obtain that
\be
&&\mE[F(X_t)|\mathscr{F}_t^Y]Z_t\no\\
&=&\int_0^t\mE[F(X_s)|\mathscr{F}_s^Y]\dif Z_s+\int_0^tZ_s\dif \mE[F(X_s)|\mathscr{F}_s^Y]+\int_0^t\int_{\mU_2}D(s,u)\(1-L(u)\)N_\lambda(\dif s, \dif u)\no\\
&=&\int_0^t\int_{\mU_2}\(1-L(u)\)\mE[F(X_s)|\mathscr{F}_s^Y]\left(\mP_{s-}\left(\lambda(s,\cdot,u)\right)-1\right)\nu_2(\dif u)\dif s\no\\
&&+\int_0^t\int_{\mU_2}D(s,u)\(1-L(u)\)\mP_{s-}\left(\lambda(s,\cdot,u)\right)\nu_2(\dif u)\dif s\no\\
&&+\int_0^tZ_s\mE[(\cL_s F)(X_s)|\mathscr{F}_s^Y]\dif s+\Sigma^4_t,
\label{disc2}
\ee
where $(\Sigma^4_t)$ is an $(\mathscr{F}^Y_t)$-adapted local martingale and given by
\ce
\Sigma^4_t&:=&\int_0^t\int_{\mU_2}\(1-L(u)\)\mE[F(X_s)|\mathscr{F}_s^Y]\tilde{\bar{N}}(\dif s, \dif u)+\int_0^t\int_{\mU_2}D(s,u)\(1-L(u)\)\tilde{\bar{N}}(\dif s, \dif u)\no\\
&&+\int_0^tZ_sE_s\dif\bar{W}_s+\int_0^t\int_{\mU_2}Z_sD(s,u)\tilde{\bar{N}}(\dif s, \dif u).
\de

Comparing (\ref{disc1}) with (\ref{disc2}), we know that
\be
D(s,u)=\frac{\mP_{s-}\(F\lambda(s,\cdot,u)\)-\mP_{s-}(F)\mP_{s-}\left(\lambda(s,\cdot,u)\right)}{\mP_{s-}\left(\lambda(s,\cdot,u)\right)}, \quad a.s. \mP,
\label{dpro}
\ee
where $\mP_{s-}\(F\lambda(s,\cdot,u)\)$ and $\mP_{s-}(F)$ are the $(\mP, \mathscr{F}_s^Y)$-predictable projections of $F(X_s)\lambda(s,X_s,u)$ and $F(X_s)$, respectively.

Finally, we attain (\ref{kseq0}) by replacing $E_s$ and $D(s,u)$ in (\ref{cdun}) with (\ref{cpro}) and (\ref{dpro}). Thus, the proof is complete.
\end{proof}

\subsection{The Zakai equation}\label{zaka}

Set
\ce
\tilde{\mP}_t(F):=\tilde{\mE}[F(X_t)\Lambda_t|\mathscr{F}_t^Y], \quad F\in\cB(\mR^n),
\de
where $\tilde{\mE}$ denotes expectation under the measure $\tilde{\mP}$. The equation satisfied by $\tilde{\mP}_t(F)$ is called the Zakai equation. In the following, we deduce the Zakai equation.

\bt\label{zakait}(The Zakai equation)\\
The Zakai equation of the system (\ref{Eq1}) is given by
\be
\tilde{\mP}_t(F)&=&\tilde{\mP}_0(F)+\int_0^t\tilde{\mP}_s(\cL_s F)\dif s
+\sum\limits_{i=1}^n\sum\limits_{l=1}^m\int_0^t\tilde{\mP}_s\left(Fh^l(s,\cdot,Y_s)+\frac{\partial F(\cdot)}{\partial x_i}\sigma^{il}_1(s,\cdot)\right)\dif \tilde{W}^l_s\no\\
&&+\int_0^t\int_{\mU_2}\tilde{\mP}_{s-}\left(F(\lambda(s,\cdot,u)-1)\right)\tilde{N}(\dif s, \dif u), \quad F\in\cC_c^\infty(\mR^n), \quad t\in[0,T].
\label{zakaieq0}
\ee
\et
\begin{proof}
Note that 
$$
\mP_t(F)=\mE[F(X_t)|\mathscr{F}_t^Y], \quad \tilde{\mP}_t(F)=\tilde{\mE}[F(X_t)\Lambda_t|\mathscr{F}_t^Y].
$$
So, by the Kallianpur-Striebel formula it holds that
\be
\mP_t(F)=\mE[F(X_t)|\mathscr{F}_t^Y]=\frac{\tilde{\mE}[F(X_t)\Lambda_t|\mathscr{F}_t^Y]}
{\tilde{\mE}[\Lambda_t|\mathscr{F}_t^Y]}=\frac{\tilde{\mP}_t(F)}{\tilde{\mP}_t(1)}.
\label{bayf}
\ee
Therefore, $\mP_t(F)\tilde{\mP}_t(1)=\tilde{\mP}_t(F)$ and then the equation which $\mP_t(F)\tilde{\mP}_t(1)$ satisfies is exactly the Zakai equation. 

First of all, we search for the equation which $\tilde{\mP}_t(1)$ satisfies. Note that
\ce
\Lambda_t&=&\exp\bigg\{\int_0^th^i(s,X_s,Y_s)\dif W^i_s+\frac{1}{2}\int_0^t
\left|h(s,X_s,Y_s)\right|^2\dif s
+\int_0^t\int_{\mU_2}\log\lambda(s,X_{s-},u)\tilde{N}_{\lambda}(\dif s, \dif u)\\
&&\quad\qquad +\int_0^t\int_{\mU_2}\(1-\lambda(s,X_s,u)-\lambda(s,X_s,u)\log\lambda(s,X_s,u)\)\nu_2(\dif u)\dif s\bigg\}.
\de
And then by the It\^o formula, one can obtain that
\ce
\Lambda_t&=&1+\int_0^t\Lambda_sh^i(s,X_s,Y_s)\dif W^i_s
+\frac{1}{2}\int_0^t\Lambda_s\left|h(s,X_s,Y_s)\right|^2\dif s+\frac{1}{2}\int_0^t\Lambda_s\left|h(s,X_s,Y_s)\right|^2\dif s\\
&&+\int_0^t\int_{\mU_2}\Lambda_{s-}\left(\lambda(s,X_{s-},u)-1\right)N_\lambda(\dif s, \dif u)
+\int_0^t\int_{\mU_2}\Lambda_{s-}(1-\lambda(s,X_{s-},u))\nu_2(\dif u)\dif s\\
&=&1+\int_0^t\Lambda_sh^i(s,X_s,Y_s)\dif \tilde{W}^i_s
+\int_0^t\int_{\mU_2}\Lambda_{s-}\left(\lambda(s,X_{s-},u)-1\right)\tilde{N}(\dif s, \dif u).
\de
Taking the conditional expectation with respect to $\mathscr{F}_t^Y$ under the measure $\tilde{\mP}$, by \cite[Theorem 1.4.7]{blr} we have that 
\ce
\tilde{\mE}[\Lambda_t|\mathscr{F}_t^Y]&=&1+\int_0^t\tilde{\mE}\left[\Lambda_sh^i(s,X_s,Y_s)|\mathscr{F}_s^Y\right]\dif \tilde{W}^i_s\\
&&+\int_0^t\int_{\mU_2}\tilde{\mE}\left[\Lambda_{s-}\left(\lambda(s,X_{s-},u)-1\right)|\mathscr{F}_s^Y\right]\tilde{N}(\dif s, \dif u),
\de
i.e.
\ce
\tilde{\mP}_t(1)=1+\int_0^t\tilde{\mP}_s(1)\mP_s(h^i(s,\cdot,Y_s))\dif \tilde{W}^i_s+\int_0^t\int_{\mU_2}\tilde{\mP}_{s-}(1)\mP_{s-}(\lambda(s,\cdot,u)-1)\tilde{N}(\dif s, \dif u).
\de

Next, applying the It\^o formula to $\mP_t(F)\tilde{\mP}_t(1)$, one can get that 
\ce
&&\mP_t(F)\tilde{\mP}_t(1)\\
&=&\mP_0(F)+\int_0^t\mP_s(F)\dif \tilde{\mP}_s(1)+\int_0^t\tilde{\mP}_s(1)\dif \mP_s(F)+\int_0^t\tilde{\mP}_s(1)\mP_s(h^i(s,\cdot,Y_s))E^i_s\dif s\\
&&+\int_0^t\int_{\mU_2}\tilde{\mP}_{s-}(1)\mP_{s-}(\lambda(s,\cdot,u)-1)D(s,u)N_{\lambda}(\dif s, \dif u)\\
&=&\mP_0(F)+\int_0^t\mP_s(F)\tilde{\mP}_s(1)\mP_s(h^i(s,\cdot,Y_s))\dif \tilde{W}^i_s+\int_0^t\tilde{\mP}_s(1)\mP_s(\cL_s F)\dif s+\int_0^t\tilde{\mP}_s(1)E^i_s\dif \tilde{W}^i_s\\
&&+\int_0^t\int_{\mU_2}\mP_{s-}(F)\tilde{\mP}_{s-}(1)\mP_{s-}(\lambda(s,\cdot,u)-1)\tilde{N}(\dif s, \dif u)-\int_0^t\tilde{\mP}_s(1)E^i_s\mP_s(h^i(s,\cdot,Y_s))\dif s\\
&&+\int_0^t\int_{\mU_2}\tilde{\mP}_s(1)D(s,u)\tilde{N}(\dif s, \dif u)-\int_0^t\int_{\mU_2}\tilde{\mP}_s(1)D(s,u)\mP_s(\lambda(s,\cdot,u)-1)\nu_2(\dif u)\dif s\\
&&+\int_0^t\tilde{\mP}_s(1)\mP_s(h^i(s,\cdot,Y_s))E^i_s\dif s+\int_0^t\int_{\mU_2}\tilde{\mP}_{s-}(1)\mP_{s-}(\lambda(s,\cdot,u)-1)D(s,u)\tilde{N}(\dif s, \dif u)\\
&&+\int_0^t\int_{\mU_2}\tilde{\mP}_s(1)\mP_s(\lambda(s,\cdot,u)-1)D(s,u)\nu_2(\dif u)\dif s\\
&=&\mP_0(F)+\int_0^t\tilde{\mP}_s(1)\mP_s(\cL_s F)\dif s
+\int_0^t\tilde{\mP}_s(1)\mP_s\left(Fh^l(s,\cdot,Y_s)+\frac{\partial F(\cdot)}{\partial x_i}\sigma^{il}_1(s,\cdot)\right)\dif \tilde{W}^l_s\\
&&+\int_0^t\int_{\mU_2}\tilde{\mP}_{s-}(1)\mP_s\left(F(\lambda(s,\cdot,u)-1)\right)\tilde{N}(\dif s, \dif u),
\de
where $E_s$ and $D(s,u)$ are given by (\ref{cpro}) and (\ref{dpro}), respectively. And then rewriting the above equation, by (\ref{bayf}) we finally obtain (\ref{zakaieq0}). The proof is over.
\end{proof}

\section{The pathwise uniqueness and uniqueness in joint law of weak solutions for the Zakai equation}\label{unzaks}

In the section we firstly define weak solutions, the pathwise uniqueness and uniqueness in joint law of weak solutions for the Zakai equation. And then we show the pathwise uniqueness and uniqueness in joint law  for weak solutions to the Zakai equation by means of a family of operators.

Let $\cP(\mR^n)$ denote the set of the probability measures on $\mR^n$ and $\cM^+(\mR^n)$ denote the set of positive bounded Borel measures on $\mR^n$. Let $\cM(\mR^n)$ denote the set of finite signed measures on $\mR^n$. For a process $\mu_{\cdot}$ valued in $\cP(\mR^n)$, $\cM^+(\mR^n)$ or $\cM(\mR^n)$, $<\mu_t, F>:=\int_{\mR^n}
F(x)\mu_t(\cdot,\dif x), F\in\cB(\mR^n)$.

Let $\mH$ be the collection of all square-integrable functions on $\mR^n$ with the norm $\|F\|^2_{\mH}=\int_{\mR^n}|F(x)|^2\dif x$ and the inner product $<F_1, F_2>_{\mH}=\int_{\mR^n}F_1(x)F_2(x)\dif x$ for $F, F_1, F_2\in\mH$. Let  $\{\phi_j, j=1,2,...\}$ be a complete orthogonal basis in $\mH$. For $\mu\in\cM(\mR^n)$, $\mu\in\mH$ means that 
$$
\|\mu\|_{\mH}^2:=\sum\limits_{j=1}^\infty|\<\mu, \phi_j\>|^2<\infty. 
$$
And then we define a family of operators on $\mH$. For $\e>0$, set
\ce
&&(S_{\e}\mu)(x):=\int_{\mR^n}(2\pi \e)^{-\frac{n}{2}}\exp\left\{-\frac{|x-y|^2}{2\e}\right\}\mu(\dif y), \quad \mu\in\cM(\mR^n),\\
&&(S_{\e}F)(x):=\int_{\mR^n}(2\pi \e)^{-\frac{n}{2}}\exp\left\{-\frac{|x-y|^2}{2\e}\right\}F(y)\dif y, \quad F\in\mH,
\de
and then one can justify $S_{\e}\mu, S_{\e}F\in\mH$. Moreover, we collect some properties of $S_{\e}$ in the following lemmas.

\bl\label{prose}
For $\mu\in\cM(\mR^n)$, $\e>0$ and $F\in\mH$,

(i) $\|S_{2\e}|\mu|\|_{\mH}\leq \|S_{\e}|\mu|\|_{\mH}$, where $|\mu|$ stands for the total variation  measure of $\mu$.

(ii) $<S_{\e}\mu, F>_{\mH}=<\mu, S_{\e}F>$.

(iii) If $\frac{\partial F}{\partial x_i}\in\mH$, 
$$
\frac{\partial (S_{\e}F)}{\partial x_i}=S_{\e}\frac{\partial F}{\partial x_i}.
$$
\el

\bl\label{onede}
Let $\xi\in\cM(\mR^n)$.

(i) Suppose that $\psi: \mR^n\rightarrow \mR$  satisfies 
$$
|\psi(x)|\leq C_1, \quad x\in\mR^n.
$$
Then there exists a positive constant $C_2$ such that 
\ce
\|S_{\e}(\psi\xi)\|_{\mH}\leq C_2\|S_{\e}(|\xi|)\|_{\mH}.
\de

(ii) Suppose that $\psi_i: \mR^n\rightarrow \mR$, $i=1,2$, satisfy 
\ce
&&|\psi_i(x_1)-\psi_i(x_2)|\leq C_3|x_1-x_2|, \quad x_1, x_2\in\mR^n,\\
&&|\psi_i(x)|\leq C_3, \quad x\in\mR^n.
\de
Then there exists a positive constant $C_4$ only depending on $\psi_1, \psi_2$ such that 
\ce
|<S_{\e}(\psi_1\xi), \frac{\partial}{\partial x_i}S_{\e}(\psi_2\xi)>|\leq C_4\|S_{\e}(|\xi|)\|^2_{\mH}.
\de
\el

\bl\label{secde}
Assume that $\Psi_1: \mR^n\rightarrow \mR^n\times\mR^d, \Psi_2: \mR^n\rightarrow \mR^n\times\mR^m$, satisfy 
\ce
&&\|\Psi_i(x_1)-\Psi_i(x_2)\|\leq C_5|x_1-x_2|, \quad x_1, x_2\in\mR^n,\\
&&\|\Psi_i(x)\|\leq C_6, \quad x\in\mR^n, \quad i=1,2,
\de
and $\xi\in\cM(\mR^n)$. Then there exist two positive constants $C_7, C_8$ depending on $\Psi_1, \Psi_2$, respectively, such that
\ce
&&<S_{\e}\xi,\frac{\partial^2}{\partial x_k\partial x_j}
S_{\e}(\Psi_1^{kl}(\cdot) \Psi_1^{lj}(\cdot)\xi)>_{\mH}+\sum_{l=1}^d\|\frac{\partial}{\partial x_i}S_{\e}(\Psi^{il}_1(\cdot)\xi)\|_{\mH}^2\leq C_7\|S_{\e}(|\xi|)\|^2_{\mH},\\
&&<S_{\e}\xi,\frac{\partial^2}{\partial x_k\partial x_j}S_{\e}(\Psi_2^{kl}(\cdot) \Psi_2^{lj}(\cdot)\xi)>_{\mH}+\sum_{l=1}^m\|\frac{\partial}{\partial x_i}S_{\e}(\Psi^{il}_2(\cdot)\xi)\|_{\mH}^2\leq C_8\|S_{\e}(|\xi|)\|^2_{\mH}.
\de
\el

Since the proofs of the above lemmas are direct, we omit them. Next, we study the Zakai equation (\ref{zakaieq0}).

\bd\label{soluzakai}
$\{(\hat{\Omega}, \hat{\mathscr{F}}, \{\hat{\mathscr{F}}_t\}_{t\in[0,T]},\hat{\mP}), (\hat{\mu}_t,
\hat{W}_t, \hat{N}(\dif t, \dif u))\}$ is called a weak solution of the Zakai equation (\ref{zakaieq0}), if the following holds:

(i) $(\hat{\Omega}, \hat{\mathscr{F}}, \{\hat{\mathscr{F}}_t\}_{t\in[0,T]},\hat{\mP})$ is a complete filtered
probability space;

(ii) $\hat{\mu}_t$ is a $\cM^+(\mR^n)$-valued $(\hat{\mathscr{F}}_t)$-adapted c\`adl\`ag process and $\hat{\mu}_0\in\cP(\mR^n)$;

(iii) $\hat{W}_t$ is a $m$-dimensional $(\hat{\mathscr{F}}_t)$-adapted Brownian motion;

(iv) $\hat{N}(\dif t, \dif u)$ is a Poisson random measure with a predictable compensator $\dif t\nu_2(\dif u)$;

(v) $(\hat{\mu}_t, \hat{W}_t, \hat{N}(\dif t, \dif u))$ satisfies the following equation
\be
<\hat{\mu}_t, F>&=&<\hat{\mu}_0, F>+\int_0^t<\hat{\mu}_s, \cL_s F>\dif s
+\sum\limits_{i=1}^n\sum\limits_{l=1}^m\int_0^t<\hat{\mu}_s, \frac{\partial F(\cdot)}{\partial x_i}\sigma^{il}_1(s,\cdot)>\dif \hat{W}^l_s\no\\
&&+\int_0^t\int_{\mU_2}<\hat{\mu}_{s-}, F(\lambda(s,\cdot,u)-1)>\tilde{\hat{N}}(\dif s, \dif u) \no\\
&&+\sum\limits_{l=1}^m\int_0^t<\hat{\mu}_s, Fh^l(s,\cdot,Y_s)>\dif \hat{W}^l_s,\quad t\in[0,T], \quad F\in \cC_c^\infty(\mR^n),
\label{zakaieq2}
\ee
where $\tilde{\hat{N}}(\dif t, \dif u):=\hat{N}(\dif t, \dif u)-\dif t\nu_2(\dif u)$.
\ed

By the deduction in Section \ref{nonfilter}, it is obvious that $\{(\Omega, \mathscr{F}, \{\mathscr{F}_t\}_{t\in[0,T]}, \tilde{\mP}),
(\tilde{\mP}_t, \tilde{W}_t, N_\lambda(\dif t, \dif u))\}$ is a weak solution of the Zakai equation (\ref{zakaieq0}).

\bd\label{paunzakai}
The pathwise uniqueness of weak solutions for the Zakai equation (\ref{zakaieq0}) means that if there exist two weak solutions $\{(\hat{\Omega}, \hat{\mathscr{F}}, \{\hat{\mathscr{F}}_t\}_{t\in[0,T]}, \hat{\mP}), (\hat{\mu}^1_t,\hat{W}_t, \hat{N}(\dif t, \dif u))\}$ and $\{(\hat{\Omega}, \hat{\mathscr{F}}, \{\hat{\mathscr{F}}_t\}_{t\in[0,T]}, \hat{\mP}), (\hat{\mu}^2_t,
\hat{W}_t, \hat{N}(\dif t, \dif u))\}$ with $\hat{\mP}\{\hat{\mu}^1_0=\hat{\mu}^2_0\}=1$, then
$$
\hat{\mu}^1_t=\hat{\mu}^2_t, \quad t\in[0,T], ~ a.s.\hat{\mP}.
$$
\ed

\bd\label{launzakai}
The uniqueness in joint law of weak solutions for the Zakai equation (\ref{zakaieq0}) means that if there exist two weak solutions $\{(\hat{\Omega}^1, \hat{\mathscr{F}}^1, \{\hat{\mathscr{F}}^1_t\}_{t\in[0,T]}, \hat{\mP}^1), (\hat{\mu}^1_t,\hat{W}^1_t, \hat{N}^1(\dif t, \dif u))\}$ and $\{(\hat{\Omega}^2, \hat{\mathscr{F}}^2, \{\hat{\mathscr{F}}^2_t\}_{t\in[0,T]}, \hat{\mP}^2),(\hat{\mu}^2_t,
\hat{W}^2_t, \hat{N}^2(\dif t, \dif u))\}$ with $\hat{\mP}^1\circ(\hat{\mu}^1_0)^{-1}=\hat{\mP}^2\circ(\hat{\mu}^2_0)^{-1}$, then
$\{(\hat{\mu}^1_t,\hat{W}^1_t, \hat{N}^1(\dif t, \dif u)), t\in[0,T]\}$ and $\{(\hat{\mu}^2_t,\hat{W}^2_t, \hat{N}^2(\dif t, \dif u)), t\in[0,T]\}$
have the same finite-dimensional distributions.
\ed

Applying these operators $S_{\e}$ to weak solutions for the Zakai equation (\ref{zakaieq0}), we get the following result.

\bl\label{trex}
Assume that $\{(\hat{\Omega}, \hat{\mathscr{F}}, \{\hat{\mathscr{F}}_t\}_{t\in[0,T]},\hat{\mP}), (\hat{\mu}_t,
\hat{W}_t, \hat{N}(\dif t, \dif u))\}$ is a weak solution for the Zakai equation (\ref{zakaieq0}). Set $Z_t^\e:=S_{\e}\hat{\mu}_t$, and then it holds that
\be
\hat{\mE}\|Z_t^\e\|^2_{\mH}&=&\|Z_0^\e\|^2_{\mH}-2\int_0^t\hat{\mE}<Z_s^\e, \frac{\partial}{\partial x_i}S_{\e}(b^i_1(s,\cdot)\hat{\mu}_s)>_{\mH}\dif s\no\\
&&+\int_0^t\hat{\mE}<Z_s^\e, \frac{\partial^2}{\partial x_i\partial x_j}S_{\e}(\sigma_0^{ik}(s,\cdot) \sigma_0^{kj}(s,\cdot)\hat{\mu}_s)>_{\mH}\dif s\no\\
&&+\int_0^t\hat{\mE}<Z_s^\e,\frac{\partial^2}{\partial x_i\partial x_j}
S_{\e}(\sigma_1^{il}(s,\cdot) \sigma_1^{lj}(s,\cdot)\hat{\mu}_s)>_{\mH}\dif s\no\\
&&+2\int_0^t\int_{\mU_1}\hat{\mE}\bigg[\sum_{j=1}^{\infty}<Z_s^\e, \phi_j>_{\mH}<Z_s^\e, \phi_j\big(\cdot+f_1(s,\cdot,u)\big)>_{\mH}-\|Z_s^\e\|^2_{\mH}\no\\
&&\qquad\qquad\qquad +<Z_s^\e,\frac{\partial}{\partial x_i}S_{\e}(f^i_1(s,\cdot,u)\hat{\mu}_s)>_{\mH}\bigg]\nu_1(\dif u)\dif s\no\\
&&+\sum_{l=1}^m\int_0^t\hat{\mE}\|\frac{\partial}{\partial x_i}S_{\e}(\sigma^{il}_1(s,\cdot)\hat{\mu}_s)\|_{\mH}^2\dif s
+\sum_{l=1}^m\int_0^t\hat{\mE}\|S_{\e}(h^l(s,\cdot,Y_s)\hat{\mu}_s)\|_{\mH}^2\dif s\no\\ 
&&-\int_0^t\hat{\mE}<\frac{\partial}{\partial x_i}S_{\e}(\sigma^{il}_1(s,\cdot)\hat{\mu}_s),S_{\e}(h^l(s,\cdot,Y_s)\hat{\mu}_s)>_{\mH}\dif s\no\\ 
&&+\int_0^t\int_{\mU_2}\hat{\mE}\|S_{\e}((\lambda(s,\cdot,u)-1)\hat{\mu}_s)\|_{\mH}^2\nu_2(\dif u)\dif s.
\label{extrop}
\ee
\el

To make the content compact, we prove Lemma \ref{trex} in the appendix. Besides, to obtain the uniqueness for weak solutions to the Zakai equation (\ref{zakaieq0}), we also need following stronger assumptions:

 \begin{enumerate}[($\mathbf{H}^{1'}_{b_1, \sigma_0, \sigma_1, f_1}$)] 
\item There exist an increasing function $L'_1(t)>0$ and a function $G_1: \mU_1\rightarrow\mR^+$ satisfying $\int_{\mU_1}G_1(u)\nu_1(\dif u)+\int_{\mU_1}G^2_1(u)\nu_1(\dif u)+\int_{\mU_1}G^4_1(u)\nu_1(\dif u)<\infty$ such that for $t\in[0,T]$ and $x_1, x_2\in\mR^n$, $u\in\mU_1$,
\ce
&|b_1(t,x_1)-b_1(t,x_2)|\leq L'_1(t)|x_1-x_2|,\\
&\|\sigma_0(t,x_1)-\sigma_0(t,x_2)\|\leq L'_1(t)|x_1-x_2|,\\
&\|\sigma_1(t,x_1)-\sigma_1(t,x_2)\|\leq L'_1(t)|x_1-x_2|,\\
&|f_1(t,x_1,u)-f_1(t,x_2,u)|\leq L'_1(t)G_1(u)|x_1-x_2|.
\de
\end{enumerate}
\begin{enumerate}[($\mathbf{H}^{2'}_{b_1, \sigma_0, \sigma_1, f_1}$)]
\item There exist an increasing function $K'_1(t)>0$ and a function $G_2: \mU_1\rightarrow\mR^+$ satisfying $\int_{\mU_1}G^2_2(u)\nu_1(\dif u)<\infty$ such that for $t\in[0,T]$ and $x\in\mR^n$, $u\in\mU_1$,
\ce
&&|b_1(t,x)|+\|\sigma_0(t,x)\|+\|\sigma_1(t,x)\|\leq K'_1(t),\\
&&|f_1(t,x,u)|\leq K'_1(t)G_2(u).
\de
\end{enumerate}
\begin{enumerate}[($\mathbf{H}^3_{f_1}$)]
\item For any $t\in[0,T]$ and $u\in\mU_1$, $x\to f_1(t,x,u)$ is invertible and differentiable, and there exists a function $G_3: \mU_1\rightarrow(0, \infty)$ satisfying $\int_{\mU_1}G_3(u)\nu_1(\dif u)<\infty$ such that 
\ce
\left|\det(J_{f_1}+I_n)\right|>\frac{1}{G_3(u)},
\de
where $J_{f_1}$ denotes the Jacobian matrix of $f_1(t,x,u)$ with respect to $x$ and $I_n$ is the $n$-order unit matrix.
\end{enumerate}
\begin{enumerate}[($\mathbf{H}^{3}_{b_2}$)] 
\item There exists an increasing function $L_3(t)>0$ such that for $t\in[0,T]$, $x_1, x_2\in\mR^n$  and $y\in\mR^m$,
\ce
|b_2(t,x_1,y)-b_2(t,x_2,y)|\leq L_3(t)|x_1-x_2|.
\de
\end{enumerate}

Next, we investigate the following moment property of weak solutions for the Zakai equation (\ref{zakaieq0}) under the above assumptions.

\bl\label{prozak}
Suppose that ($\mathbf{H}^{1'}_{b_1, \sigma_0, \sigma_1, f_1}$) ($\mathbf{H}^{2'}_{b_1, \sigma_0, \sigma_1, f_1}$) ($\mathbf{H}^3_{f_1}$) ($\mathbf{H}^1_{\sigma_2, f_2}$) ($\mathbf{H}^2_{b_2, \sigma_2, f_2}$) ($\mathbf{H}^{3}_{b_2}$) hold. Then for a weak solution $\{(\hat{\Omega}, \hat{\mathscr{F}}, \{\hat{\mathscr{F}}_t\}_{t\in[0,T]},\hat{\mP}), (\hat{\mu}_t,
\hat{W}_t, \hat{N}(\dif t, \dif u))\}$ of the Zakai equation (\ref{zakaieq0}) with $\hat{\mu}_0\in\mH$, it holds that for $Z_t^\e:=S_{\e}\hat{\mu}_t$,
\be\label{moin}
\hat{\mE}\|Z_t^\e\|^2_{\mH}\leq\|Z_0^\e\|^2_{\mH}+C\int_0^t\hat{\mE}\|Z^\e_s\|^2_{\mH}\dif s,
\ee
and $\hat{\mu}_t\in\mH, a.s.$ for $t\in[0,T]$.
\el

Since the proof of the above lemma is too long, we place it in the appendix. Now, it is the position to state and prove the result on the pathwise uniqueness for weak solutions to the Zakai equation. We recall that $\{(\Omega, \mathscr{F}, \{\mathscr{F}_t\}_{t\in[0,T]}, \tilde{\mP}), (\tilde{\mP}_t, \tilde{W}_t, N_\lambda(\dif t, \dif u))\}$ is a weak solution of the Zakai equation (\ref{zakaieq0}).

\bt\label{unizak}(The pathwise uniqueness)\\
Suppose that ($\mathbf{H}^{1'}_{b_1, \sigma_0, \sigma_1, f_1}$) ($\mathbf{H}^{2'}_{b_1, \sigma_0, \sigma_1, f_1}$) ($\mathbf{H}^3_{f_1}$) ($\mathbf{H}^1_{\sigma_2, f_2}$) ($\mathbf{H}^2_{b_2, \sigma_2, f_2}$) ($\mathbf{H}^{3}_{b_2}$) hold. If $\{(\Omega, \mathscr{F}, \\\{\mathscr{F}_t\}_{t\in[0,T]}, \tilde{\mP}),
(\tilde{\mu}_t, \tilde{W}_t, N_\lambda(\dif t, \dif u))\}$ with $\tilde{\mu}_0=\tilde{\mP}_0$ is another weak solution for the Zakai equation (\ref{zakaieq0}), then $\tilde{\mu}_t=\tilde{\mP}_t$ for any $t\in[0,T]$ a.s. $\tilde{\mP}$.
\et
\begin{proof}
Set $D_t:=\tilde{\mu}_t-\tilde{\mP}_t$, and then $D_t$ satisfies Eq.(\ref{zakaieq2}) due to linearity of the Zakai equation. So, by the same deduction to that in Lemma \ref{prozak}, it holds that
\ce
\tilde{\mE}\|S_\e D_t\|^2_{\mH}\leq C\int_0^t\tilde{\mE}\|S_\e(|D_s|)\|^2_{\mH}\dif s\leq C\int_0^t\tilde{\mE}\||D_s|\|^2_{\mH}\dif s=C\int_0^t\tilde{\mE}\|D_s\|^2_{\mH}\dif s.
\de
As $\e\rightarrow0$, we have that
$$
\tilde{\mE}\|D_t\|^2_{\mH}\leq C\int_0^t\tilde{\mE}\|D_s\|^2_{\mH}\dif s.
$$
Thus, it follows from the Gronwall inequality that $\tilde{\mu}_t=\tilde{\mP}_t$, a.s. for any $t\in[0,T]$. Thus, the c\`adl\`ag property of
$\tilde{\mu}_t, \tilde{\mP}_t$ in $t$ admits us to get the pathwise uniqueness, i.e.
$$
\tilde{\mu}_t=\tilde{\mP}_t, ~ \forall t\in[0,T], ~ a.s.\tilde{\mP}.
$$
The proof is complete.
\end{proof}

\bt\label{unilawzak}(The uniqueness in joint law)\\
Assume that ($\mathbf{H}^{1'}_{b_1, \sigma_0, \sigma_1, f_1}$) ($\mathbf{H}^{2'}_{b_1, \sigma_0, \sigma_1, f_1}$) ($\mathbf{H}^3_{f_1}$) ($\mathbf{H}^1_{\sigma_2, f_2}$) ($\mathbf{H}^2_{b_2, \sigma_2, f_2}$) ($\mathbf{H}^{3}_{b_2}$) hold. Then weak solutions of the Zakai equation (\ref{zakaieq0}) have the uniqueness in joint law.
\et

Since the proof of the above theorem is similar to that of Theorem 4 (ii)  in \cite{q2}, we omit it.

\section{The uniqueness in joint law of weak solutions for the Kushner-Stratonovich equation}\label{unks}

In the section, we introduce weak solutions and the uniqueness in joint law of weak solutions for the Kushner-Stratonovich equation (\ref{kseq0}). And then, the uniqueness in joint law for the Kushner-Stratonovich equation (\ref{kseq0}) is proved through the relationship between weak solutions of the Zakai equation and that of the Kushner-Stratonovich equation.

\bd\label{soluks}
If there exists the pair $\{(\bar{\Omega}, \bar{\mathscr{F}}, \{\bar{\mathscr{F}}_t\}_{t\in[0,T]}, \bar{\mP}), (\pi_t,
I_t, U(\dif t, \dif u))\}$ such that the following holds:

(i) $(\bar{\Omega}, \bar{\mathscr{F}}, \{\bar{\mathscr{F}}_t\}_{t\in[0,T]},\bar{\mP})$ is a complete filtered
probability space;

(ii) $\pi_t$ is a $\cP(\mR^n)$-valued $(\bar{\mathscr{F}}_t)$-adapted c\`adl\`ag process;

(iii) $I_t$ is a $m$-dimensional $(\bar{\mathscr{F}}_t)$-adapted Brownian motion;

(iv) $U(\dif t, \dif u)$ is a Poisson random measure with a predictable compensator \\
$\pi_{t}\left(\lambda(t,\cdot,u)\right)\dif t\nu_2(\dif u)$;

(v) $(\pi_t, I_t, U(\dif t, \dif u))$ satisfies the following equation
\be
<\pi_t, F>&=&<\pi_0, F>+\int_0^t<\pi_s, \cL_s F>\dif s+\sum\limits_{i=1}^n\sum\limits_{l=1}^m\int_0^t<\pi_s, \frac{\partial F(\cdot)}{\partial x_i}\sigma^{il}_1(s,\cdot)>\dif I^l_s\no\\
&&+\sum\limits_{l=1}^m\int_0^t\bigg(<\pi_s, Fh^l(s,\cdot,Y_s)>-<\pi_s, F><\pi_s, h^l(s,\cdot,Y_s)>\bigg)\dif I^l_s\no\\
&&+\int_0^t\int_{\mU_2}\frac{<\pi_{s-}, F\lambda(s,\cdot,u)>-<\pi_{s-}, F><\pi_{s-}
,\lambda(s,\cdot,u)>}{<\pi_{s-}, \lambda(s,\cdot,u)>}\tilde{U}(\dif s, \dif u), \no\\
&&\qquad\qquad\qquad\qquad t\in[0,T], \quad F\in\cC_c^\infty(\mR^n),
\label{kseq}
\ee
where
$$
\tilde{U}(\dif t, \dif u)=U(\dif t, \dif u)-\pi_{t}\left(\lambda(t,\cdot,u)\right)\dif t\nu_2(\dif u),
$$
then $\{(\bar{\Omega}, \bar{\mathscr{F}}, \{\bar{\mathscr{F}}_t\}_{t\in[0,T]}, \bar{\mP}), (\pi_t,
I_t, U(\dif t, \dif u))\}$ is called a weak solution of the Kushner-Stratonovich equation (\ref{kseq0}).
\ed

By the deduction in Section \ref{nonfilter}, it is obvious that $\{(\Omega, \mathscr{F}, \{\mathscr{F}_t\}_{t\in[0,T]}, \mP),
(\mP_t, \bar{W}_t, N_\lambda(\dif t, \dif u))\}$ is a weak solution of the Kushner-Stratonovich equation (\ref{kseq0}).

\bd\label{launks}
The uniqueness in joint law of weak solutions for the Kushner-Stratonovich equation (\ref{kseq0}) means that if there exist two weak solutions $\{(\bar{\Omega}^1, \bar{\mathscr{F}}^1, \{\bar{\mathscr{F}}^1_t\}_{t\in[0,T]}, \bar{\mP}^1), (\pi^1_t, I^1_t,\\ U^1(\dif t, \dif u))\}$ and $\{(\bar{\Omega}^2, \bar{\mathscr{F}}^2, \{\bar{\mathscr{F}}^2_t\}_{t\in[0,T]}, \bar{\mP}^2), (\pi^2_t, I^2_t, U^2(\dif t, \dif u))\}$ with $\bar{\mP}^1\circ(\pi^1_0)^{-1}=\bar{\mP}^2\circ(\pi^2_0)^{-1}$, then 
$\{(\pi^1_t, I^1_t, U^1(\dif t, \dif u)), t\in[0,T]\}$ and $\{(\pi^2_t, I^2_t, U^2(\dif t, \dif u)), t\in[0,T]\}$
have the same finite-dimensional distributions.
\ed

Here, we give out the main result in the section.

\bt\label{uniks}(The uniqueness in joint law)\\
Suppose that ($\mathbf{H}^{1'}_{b_1, \sigma_0, \sigma_1, f_1}$) ($\mathbf{H}^{2'}_{b_1, \sigma_0, \sigma_1, f_1}$) ($\mathbf{H}^3_{f_1}$) ($\mathbf{H}^1_{\sigma_2, f_2}$) ($\mathbf{H}^2_{b_2, \sigma_2, f_2}$) ($\mathbf{H}^{3}_{b_2}$) hold. Then weak solutions of the Kushner-Stratonovich equation (\ref{kseq0}) have the uniqueness in joint law.
\et

Since the proof of the above theorem is similar to that of Theorem 5  in \cite{q2}, we omit it.

\br\label{notun}
Since in Definition \ref{soluks} $U(\dif t, \dif u)$ depends on $\pi_t$, usual strong solutions and the pathwise uniqueness of strong solutions for the Kushner-Stratonovich equation can not be defined. Thus, we don't consider its pathwise uniqueness here. 
\er

\section{Conclusion}\label{con}

In the paper, we consider nonlinear filtering problems of stochastic differential equations driven by correlated L\'evy noises. First, we establish the Kushner-Stratonovich and Zakai equations by martingale problems and the Kallianpur-Striebel formula. And then, the pathwise uniqueness and uniqueness in joint law of weak solutions for the Zakai equation are shown. Finally, we study the uniqueness in joint law of weak solutions for the Kushner-Stratonovich equation.

Our method also can be used to solve nonlinear filtering problems of stochastic differential equations with correlated sensor L\'evy noises. Concretely speaking, fix $T>0$ and consider the following signal-observation system $(\check{X}_t, \check{Y}_t)$ on $\mR^n\times\mR^m$:
\be\left\{\begin{array}{l}
\dif \check{X}_t=\check{b}_1(t,\check{X}_t)\dif t+\check{\sigma}_1(t,\check{X}_t)\dif W_t+\int_{\mU_1}\check{f}_1(t,\check{X}_{t-},u)\tilde{N_p}(\dif t, \dif u),\\
\dif \check{Y}_t=\check{b}_2(t,\check{X}_t,\check{Y}_t)\dif t+\check{\sigma}_2\dif W_t+\check{\sigma}_3\dif B_t+\int_{\mU_2}\check{f}_2(t, \check{Y}_{t-}, u)\tilde{\check{N}}_{\lambda}(\dif t, \dif u), 
 \end{array}
\right. \no\\
0\leq t\leq T,
\label{Eq0} 
\ee
where $\check{N}_{\lambda}(\dif t,\dif u)$ is an integer-valued random measure, its predictable compensator is given by $\lambda(t,\check{X}_{t-},u)\dif t\nu_2(\dif u)$ and $\tilde{\check{N}}_{\lambda}(\dif t, \dif u):=\check{N}_{\lambda}(\dif t,\dif u)-\lambda(t,\check{X}_{t-},u)\dif t\nu_2(\dif u)$. The initial value $\check{X}_0$ is assumed to be a square integrable random variable independent of $\check{Y}_0, W_{\cdot}, B_{\cdot}, N_p, \check{N}_{\lambda}$. Moreover, $W_{\cdot}, B_{\cdot}, N_p, \check{N}_{\lambda}$ are mutually independent.

The mappings $\check{b}_1:[0,T]\times\mR^n\mapsto\mR^n$, $\check{\sigma}_1:[0,T]\times\mR^n\mapsto\mR^{n\times m}$, $\check{f}_1:[0,T]\times\mR^n\times\mU_1\mapsto\mR^n$, $\check{b}_2:[0,T]\times\mR^n\times\mR^m\mapsto\mR^m$ and $\check{f}_2:[0,T]\times\mR^m\times\mU_2\mapsto\mR^m$ are all Borel measurable. $\check{\sigma}_2, \check{\sigma}_3$ are $m\times m$ and $m\times d$ real matrices, respectively. Moreover, we make the following hypotheses:
\begin{enumerate}[(i)]
\item $\check{b}_1, \check{\sigma}_1, \check{f}_1$ satisfy ($\mathbf{H}^1_{b_1, \sigma_0, \sigma_1, f_1}$)-($\mathbf{H}^2_{b_1,\sigma_0,\sigma_1, f_1}$), where $\check{b}_1, \check{\sigma}_1, \check{f}_1$ replace $b_1, \sigma_1, f_1$;
\item $\check{b}_2(t,x,y)$ is bounded for all $t\in[0,T], x\in\mR^n, y\in\mR^m$;
\item $\check{\sigma}_2\check{\sigma}^*_2+\check{\sigma}_3\check{\sigma}^*_3=I_m,$ where $\check{\sigma}^*_2$ stands for the transpose of the matrix $\check{\sigma}_2$ and $I_m$ is the $m$-order unit matrix;
\item $\check{f}_2$ satisfies ($\mathbf{H}^1_{\sigma_2, f_2}$)-($\mathbf{H}^2_{b_2,\sigma_2, f_2}$), where $\check{f}_2$ replaces $f_2$.
\end{enumerate}

Under the above assumptions and ($\mathbf{H}_{\lambda}$), by the similar method to that in the proofs of Theorem \ref{ks} and \ref{zakait}, we can establish the Kushner-Stratonovich equation and Zakai equation of the system (\ref{Eq0}). Specially, a unique strong solution of the system (\ref{Eq0}) is denoted as $(\check{X}_t,\check{Y}_t)$. And then set
\ce
\check{\mP}_t(F):=\mE[F(\check{X}_t)|\mathscr{F}_t^{\check{Y}}], \quad F\in\cB(\mR^n).
\de

\bc (The Kushner-Stratonovich equation)\\  
For $F\in\cC_c^\infty(\mR^n)$, the Kushner-Stratonovich equation of the system (\ref{Eq0}) is given by
\be
\check{\mP}_t(F)&=&\check{\mP}_0(F)+\int_0^t\check{\mP}_s(\check{\cL}_s F)\dif s+\int_0^t\check{\mP}_s\left(\frac{\partial F(\cdot)}{\partial x_i}\check{\sigma}^{ik}_1(s,\cdot)\check{\sigma}^{jk}_2\right)\dif  \bar{V}^j_s\no\\
&&+\sum\limits_{l=1}^m\int_0^t\bigg(\check{\mP}_s\left(F\check{b}_2^l(s,\cdot,Y_s)\right)-\check{\mP}_s\left(F\right)\check{\mP}_s\left(\check{b}_2^l(s,\cdot,Y_s)\right)\bigg)\dif \bar{V}^l_s\no\\
&&+\int_0^t\int_{\mU_2}\frac{\check{\mP}_{s-}\left(F\lambda(s,\cdot,u)\right)-\check{\mP}_{s-}\left(F\right)\check{\mP}_{s-}
\left(\lambda(s,\cdot,u)\right)}{\check{\mP}_{s-}\left(\lambda(s,\cdot,u)\right)}\bar{\check{N}}(\dif s, \dif u), t\in[0,T],\no\\
\label{kseq01}
\ee
where 
\ce
(\check{\cL_s} F)(x)&:=&\frac{\partial F(x)}{\partial x_i}\check{b}^i_1(s,x)+\frac{1}{2}\frac{\partial^2F(x)}{\partial x_i\partial x_j}
\check{\sigma}_1^{il}(s,x) \check{\sigma}_1^{lj}(s,x)\\
&&+\int_{\mU_1}\Big[F\big(x+\check{f}_1(s,x,u)\big)-F(x)
-\frac{\partial F(x)}{\partial x_i}\check{f}^i_1(s,x,u)\Big]\nu_1(\dif u),
\de
and $\bar{V}_t:=\check{\sigma}_2 W_t+\check{\sigma}_3 B_t+\int_0^t \check{b}_2(s,\check{X}_s,\check{Y}_s)\dif s-\int_0^t\check{\mP}_s(\check{b}_2(s,\cdot,\check{Y}_s))\dif s$, $\bar{\check{N}}(\dif t, \dif u)=\check{N}_\lambda(\dif t, \dif u)-\check{\mP}_{t-}\left(\lambda(t,\cdot,u)\right)\nu_2(\dif u)\dif t$.
\ec

Next, set
\ce
&&V_t:=\check{\sigma}_2 W_t+\check{\sigma}_3 B_t, \\
&&\gamma^{-1}_t:=\exp\bigg\{-\int_0^t\check{b}_2^i(s,\check{X}_s,\check{Y}_s)\dif V^i_s-\frac{1}{2}\int_0^t
\left|\check{b}_2(s,\check{X}_s,\check{Y}_s)\right|^2\dif s\\
&&\quad\qquad\quad\qquad-\int_0^t\int_{\mU_2}\log\lambda(s,\check{X}_{s-},u)\tilde{\check{N}}_{\lambda}(\dif s, \dif u)\\
&&\quad\qquad\quad\qquad -\int_0^t\int_{\mU_2}\(1-\lambda(s,\check{X}_s,u)+\lambda(s,\check{X}_s,u)\log\lambda(s,\check{X}_{s},u)\)\nu_2(\dif u)\dif s\bigg\},
\de
and then by the similar deduction to \cite{qd}, we know that $\gamma^{-1}_t$ is an exponential martingale. In addition, define the probability measure 
$$
\frac{\dif \tilde{\check{\mP}}}{\dif \mP}:=\gamma^{-1}_T,
$$
and set
\ce
\tilde{\check{\mP}}_t(F):=\tilde{\check{\mE}}[F(\check{X}_t)\gamma_t|\mathscr{F}_t^{\check{Y}}],
\de
where $\tilde{\check{\mE}}$ stands for the expectation under the probability measure $\tilde{\check{\mP}}$.

\bc (The Zakai equation)\\
 The Zakai equation of the system (\ref{Eq0}) is given by
\be
\tilde{\check{\mP}}_t(F)&=&\tilde{\check{\mP}}_0(F)+\int_0^t\tilde{\check{\mP}}_s(\check{\cL}_s F)\dif s
+\sum\limits_{l=1}^m\int_0^t\tilde{\check{\mP}}_s\left(F\check{b}_2^l(s,\cdot,Y_s)+\frac{\partial F(\cdot)}{\partial x_i}\check{\sigma}^{ik}_1(s,\cdot)\check{\sigma}^{lk}_2\right)\dif \tilde{V}^l_s\no\\
&&+\int_0^t\int_{\mU_2}\tilde{\check{\mP}}_{s-}\left(F(\lambda(s,\cdot,u)-1)\right)\tilde{\check{N}}(\dif s, \dif u), \quad F\in\cC_c^\infty(\mR^n), \quad t\in[0,T],
\label{zakaieq01}
\ee
where $\tilde{V}_t:=V_t+\int_0^t \check{b}_2(s,\check{X}_s,\check{Y}_s)\dif s$ and $\tilde{\check{N}}(\dif t, \dif u):=\check{N}_\lambda(\dif t, \dif u)-\dif t\nu_2(\dif u)$.
\ec

Of course, we can study the pathwise uniqueness and uniqueness in joint law of weak solutions for Eq.(\ref{zakaieq01}) and Eq.(\ref{kseq01}) by the same means to that in Theorem \ref{unizak}, \ref{unilawzak}, \ref{uniks}.

\section{Appendix}\label{app}

In the section, we give out the proofs of Lemma \ref{brmoposs}, Lemma \ref{trex} and Lemma \ref{prozak}.

{\bf The proof of Lemma \ref{brmoposs}.}

By the similar proof to that in \cite[Page 323, Theorem 8.4]{roge}, we know that $(\bar{W}_t)$ is an $(\mathscr{F}^Y_t)$-adapted Brownian motion. Therefore, it is only necessary to prove that $\(\tilde{\bar{N}}(\dif t, \dif u)\)$ is an $(\mathscr{F}^Y_t)$-adapted martingale measure. 

First of all, by \cite[Proposition 3.2]{ceci} and \cite[Proposition 2.2]{cckc1}, it holds that $\mP_{t-}\left(\lambda(t,\cdot,u)\right)\nu_2(\dif u)\dif t$ is the $(\mP, \mathscr{F}_t^Y)$-predictable projection of $N_\lambda(\dif t, \dif u)$. And then we show that
$$
\mE[\tilde{\bar{N}}((0,t]\times A)|\mathscr{F}_v^Y]=\tilde{\bar{N}}((0,v]\times A), \quad 0<v<t, \quad A\in\mathscr{U}.
$$

We begin with the left side of the above equality. By the expression of $\tilde{\bar{N}}(\dif t, \dif u)$, it holds that
\be
&&\mE[\tilde{\bar{N}}((0,t]\times A)|\mathscr{F}_v^Y]\no\\
&=&\mE\[N_\lambda((0,t]\times A)-\int_0^t\int_A\mP_{s-}\left(\lambda(s,\cdot,u)\right)\nu_2(\dif u)\dif s|\mathscr{F}_v^Y\]\no\\
&=&\mE\[N_\lambda((0,t]\times A)-\int_0^t\int_A\lambda(s,X_s,u)\nu_2(\dif u)\dif s|\mathscr{F}_v^Y\]\no\\
&&+\mE\[\int_0^t\int_A\lambda(s,X_s,u)\nu_2(\dif u)\dif s-\int_0^t\int_A\mP_{s-}\left(\lambda(s,\cdot,u)\right)\nu_2(\dif u)\dif s|\mathscr{F}_v^Y\]\no\\
&=:&I_1+I_2.
\label{left}
\ee
For $I_1$, note that $\(N_\lambda((0,t]\times A)-\int_0^t\int_A\lambda(s,X_s,u)\nu_2(\dif u)\dif s\)$ is an $(\mathscr{F}_t)$-adapted martingale. So, it follows from the tower property of the conditional expectation that
\be
I_1&=&\mE\[\mE\[N_\lambda((0,t]\times A)-\int_0^t\int_A\lambda(s,X_s,u)\nu_2(\dif u)\dif s|\mathscr{F}_v\]|\mathscr{F}_v^Y\]\no\\
&=&\mE\[N_\lambda((0,v]\times A)-\int_0^v\int_A\lambda(s,X_s,u)\nu_2(\dif u)\dif s|\mathscr{F}_v^Y\]\no\\
&=&N_\lambda((0,v]\times A)-\int_0^v\int_A\mE[\lambda(s,X_s,u)|\mathscr{F}_v^Y]\nu_2(\dif u)\dif s,
\label{right1}
\ee
where the measurablity of $N_\lambda((0,v]\times A)$ with respect to $\mathscr{F}_v^Y$ is used in the last equality. For $I_2$, again by the tower property of the conditional expectation we have that
\be
I_2&=&\int_0^t\int_A\mE[\lambda(s,X_s,u)|\mathscr{F}_v^Y]\nu_2(\dif u)\dif s-\int_0^t\int_A\mE[\mP_{s-}\left(\lambda(s,\cdot,u)\right)|\mathscr{F}_v^Y]\nu_2(\dif u)\dif s\no\\
&=&\int_0^v\int_A\mE[\lambda(s,X_s,u)|\mathscr{F}_v^Y]\nu_2(\dif u)\dif s+\int_v^t\int_A\mE[\lambda(s,X_s,u)|\mathscr{F}_v^Y]\nu_2(\dif u)\dif s\no\\
&&-\int_0^v\int_A\mE[\mP_{s-}\left(\lambda(s,\cdot,u)\right)|\mathscr{F}_v^Y]\nu_2(\dif u)\dif s
-\int_v^t\int_A\mE[\mP_{s-}\left(\lambda(s,\cdot,u)\right)|\mathscr{F}_v^Y]\nu_2(\dif u)\dif s\no\\
&=&\int_0^v\int_A\mE[\lambda(s,X_s,u)|\mathscr{F}_v^Y]\nu_2(\dif u)\dif s-\int_0^v\int_A\mP_{s-}\left(\lambda(s,\cdot,u)\right)\nu_2(\dif u)\dif s.
\label{right2}
 \ee
Combining (\ref{right1}) (\ref{right2}) with (\ref{left}), one can obtain 
\ce
\mE[\tilde{\bar{N}}((0,t]\times A)|\mathscr{F}_v^Y]=N_\lambda((0,v]\times A)-\int_0^v\int_A\mP_{s-}\left(\lambda(s,\cdot,u)\right)\nu_2(\dif u)\dif s=\tilde{\bar{N}}((0,v]\times A).
\de
The proof is complete.

{\bf The proof of Lemma \ref{trex}.}

{\bf Step 1.} 
We establish the following equation
\be
&&<Z_t^\e, F>_{\mH}\no\\
&=&<Z_0^\e, F>_{\mH}-\int_0^t<\frac{\partial}{\partial x_i}S_{\e}(b^i_1(s,\cdot)\hat{\mu}_s),  F>_{\mH}\dif s\no\\
&&+\frac{1}{2}\int_0^t<\frac{\partial^2}{\partial x_i\partial x_j}S_{\e}(\sigma_0^{ik}(s,\cdot) \sigma_0^{kj}(s,\cdot)\hat{\mu}_s),  F>_{\mH}\dif s\no\\
&&+\frac{1}{2}\int_0^t<\frac{\partial^2}{\partial x_i\partial x_j}
S_{\e}(\sigma_1^{il}(s,\cdot) \sigma_1^{lj}(s,\cdot)\hat{\mu}_s), F>_{\mH}\dif s\no\\
&&+\int_0^t\int_{\mU_1}\bigg[<S_{\e}\hat{\mu}_s, F\big(\cdot+f_1(s,\cdot,u)\big)>_{\mH}-<S_{\e}\hat{\mu}_s, F>_{\mH}\no\\
&&\qquad\qquad +<\frac{\partial}{\partial x_i}S_{\e}(f^i_1(s,\cdot,u)\hat{\mu}_s),  F>_{\mH}\bigg]\nu_1(\dif u)\dif s\no\\
&&-\int_0^t<\frac{\partial}{\partial x_i}S_{\e}(\sigma^{il}_1(s,\cdot)\hat{\mu}_s), F>_{\mH}\dif \hat{W}^l_s+\int_0^t<S_{\e}(h^l(s,\cdot,Y_s)\hat{\mu}_s), F>_{\mH}\dif \hat{W}^l_s\no\\
&&+\int_0^t\int_{\mU_2}<S_{\e}((\lambda(s,\cdot,u)-1)\hat{\mu}_{s-}), F>_{\mH}\tilde{\hat{N}}(\dif s, \dif u), \no\\ 
&&\quad t\in[0,T], \quad F\in\cC_c^\infty(\mR^n).
\label{esexo}
\ee

By Definition \ref{soluzakai}, we know that for $F\in\cC_c^\infty(\mR^n)$
\ce
<\hat{\mu}_t, F>&=&<\hat{\mu}_0, F>+\int_0^t<\hat{\mu}_s, \cL_s F>\dif s
+\int_0^t<\hat{\mu}_s, \frac{\partial F(\cdot)}{\partial x_i}\sigma^{il}_1(s,\cdot)>\dif \hat{W}^l_s\no\\
&&+\int_0^t\int_{\mU_2}<\hat{\mu}_{s-}, F(\lambda(s,\cdot,u)-1)>\tilde{\hat{N}}(\dif s, \dif u)\no\\
&&+\int_0^t<\hat{\mu}_s, Fh^l(s,\cdot,Y_s)>\dif \hat{W}^l_s,  \quad t\in[0,T].
\de
Replacing $F$ by $S_{\e}F$ and using Lemma \ref{prose}, we obtain that
\be
<S_{\e}\hat{\mu}_t, F>_{\mH}&=&<S_{\e}\hat{\mu}_0, F>_{\mH}
+\int_0^t<\hat{\mu}_s, \cL_s(S_{\e}F)>\dif s\no\\
&&+\int_0^t<\hat{\mu}_s, \frac{\partial (S_{\e}F)(\cdot)}{\partial x_i}\sigma^{il}_1(s,\cdot)>\dif \hat{W}^l_s\no\\
&&+\int_0^t<\hat{\mu}_s, (S_{\e}F)h^l(s,\cdot,Y_s)>\dif \hat{W}^l_s\no\\
&&+\int_0^t\int_{\mU_2}<\hat{\mu}_{s-}, (S_{\e}F)(\lambda(s,\cdot,u)-1)>\tilde{\hat{N}}(\dif s, \dif u).
\label{esex}
\ee
Note that $Z_t^\e=S_{\e}\hat{\mu}_t$. Thus, we define
\ce
&&I_1:=<\hat{\mu}_s, \cL_s(S_{\e}F)>, \qquad\qquad I_2:=<\hat{\mu}_s, \frac{\partial (S_{\e}F)(\cdot)}{\partial x_i}\sigma^{il}_1(s,\cdot)>, \\
&&I_3:=<\hat{\mu}_s, (S_{\e}F)h^l(s,\cdot,Y_s)>, \quad I_4:=<\hat{\mu}_{s-}, (S_{\e}F)(\lambda(s,\cdot,u)-1)>,
\de
and rewrite them to obtain (\ref{esexo}). 

For $I_1$, by the definition of $\cL_s$ and Lemma \ref{prose}, it holds that
\be
I_1&=&<\hat{\mu}_s, \frac{\partial (S_{\e}F)(\cdot)}{\partial x_i}b^i_1(s,\cdot)>+\frac{1}{2}<\hat{\mu}_s, \frac{\partial^2(S_{\e}F)(\cdot)}{\partial x_i\partial x_j}
\sigma_0^{ik}(s,\cdot) \sigma_0^{kj}(s,\cdot)>\no\\
&&+\frac{1}{2}<\hat{\mu}_s, \frac{\partial^2(S_{\e}F)(\cdot)}{\partial x_i\partial x_j}
\sigma_1^{il}(s,\cdot) \sigma_1^{lj}(s,\cdot)>\no\\
&&+\int_{\mU_1}\bigg[<\hat{\mu}_s, (S_{\e}F)\big(\cdot+f_1(s,\cdot,u)\big)>-<\hat{\mu}_s, S_{\e}F>\no\\
&&\qquad\qquad\qquad -<\hat{\mu}_s, \frac{\partial (S_{\e}F)(\cdot)}{\partial x_i}f^i_1(s,\cdot,u)>\bigg]\nu_1(\dif u)\no\\
&=&<b^i_1(s,\cdot)\hat{\mu}_s, \frac{\partial (S_{\e}F)(\cdot)}{\partial x_i}>+\frac{1}{2}<\sigma_0^{ik}(s,\cdot) \sigma_0^{kj}(s,\cdot)\hat{\mu}_s, \frac{\partial^2(S_{\e}F)(\cdot)}{\partial x_i\partial x_j}
>\no\\
&&+\frac{1}{2}<\sigma_1^{il}(s,\cdot) \sigma_1^{lj}(s,\cdot)\hat{\mu}_s, \frac{\partial^2(S_{\e}F)(\cdot)}{\partial x_i\partial x_j}
>\no\\
&&+\int_{\mU_1}\bigg[<\hat{\mu}_s, (S_{\e}F)\big(\cdot+f_1(s,\cdot,u)\big)>-<\hat{\mu}_s, S_{\e}F>\no\\
&&\qquad\qquad\qquad -<f^i_1(s,\cdot,u)\hat{\mu}_s, \frac{\partial (S_{\e}F)(\cdot)}{\partial x_i}>\bigg]\nu_1(\dif u)\no\\
&=&<b^i_1(s,\cdot)\hat{\mu}_s, S_{\e}\frac{\partial F(\cdot)}{\partial x_i}>+\frac{1}{2}<\sigma_0^{ik}(s,\cdot) \sigma_0^{kj}(s,\cdot)\hat{\mu}_s, S_{\e}\frac{\partial^2 F(\cdot)}{\partial x_i\partial x_j}
>\no\\
&&+\frac{1}{2}<\sigma_1^{il}(s,\cdot) \sigma_1^{lj}(s,\cdot)\hat{\mu}_s, S_{\e}\frac{\partial^2F(\cdot)}{\partial x_i\partial x_j}
>\no\\
&&+\int_{\mU_1}\bigg[<\hat{\mu}_s, (S_{\e}F)\big(\cdot+f_1(s,\cdot,u)\big)>-<\hat{\mu}_s, S_{\e}F>\no\\
&&\qquad\qquad\qquad -<f^i_1(s,\cdot,u)\hat{\mu}_s, S_{\e}\frac{\partial F(\cdot)}{\partial x_i}>\bigg]\nu_1(\dif u)\no\\
&=&<S_{\e}(b^i_1(s,\cdot)\hat{\mu}_s), \frac{\partial F(\cdot)}{\partial x_i}>_{\mH}+\frac{1}{2}<S_{\e}(\sigma_0^{ik}(s,\cdot) \sigma_0^{kj}(s,\cdot)\hat{\mu}_s), \frac{\partial^2 F(\cdot)}{\partial x_i\partial x_j}
>_{\mH}\no\\
&&+\frac{1}{2}<S_{\e}(\sigma_1^{il}(s,\cdot) \sigma_1^{lj}(s,\cdot)\hat{\mu}_s), \frac{\partial^2F(\cdot)}{\partial x_i\partial x_j}
>_{\mH}\no\\
&&+\int_{\mU_1}\bigg[<S_{\e}\hat{\mu}_s, F\big(\cdot+f_1(s,\cdot,u)\big)>_{\mH}-<S_{\e}\hat{\mu}_s, F>_{\mH}\no\\
&&\qquad\qquad\qquad -<S_{\e}(f^i_1(s,\cdot,u)\hat{\mu}_s), \frac{\partial F(\cdot)}{\partial x_i}>_{\mH}\bigg]\nu_1(\dif u)\no\\
&=&-<\frac{\partial}{\partial x_i}S_{\e}(b^i_1(s,\cdot)\hat{\mu}_s),  F>_{\mH}+\frac{1}{2}<\frac{\partial^2}{\partial x_i\partial x_j}S_{\e}(\sigma_0^{ik}(s,\cdot) \sigma_0^{kj}(s,\cdot)\hat{\mu}_s),  F>_{\mH}\no\\
&&+\frac{1}{2}<\frac{\partial^2}{\partial x_i\partial x_j}
S_{\e}(\sigma_1^{il}(s,\cdot) \sigma_1^{lj}(s,\cdot)\hat{\mu}_s), F>_{\mH}\no\\
&&+\int_{\mU_1}\bigg[<S_{\e}\hat{\mu}_s, F\big(\cdot+f_1(s,\cdot,u)\big)>_{\mH}-<S_{\e}\hat{\mu}_s, F>_{\mH}\no\\
&&\qquad\qquad\qquad +<\frac{\partial}{\partial x_i}S_{\e}(f^i_1(s,\cdot,u)\hat{\mu}_s),  F>_{\mH}\bigg]\nu_1(\dif u),
\label{i1es}
\ee
where in the last equality the formula for integration by parts is used. 

For $I_2$, it follows from Lemma \ref{prose} that 
\be
I_2&=&<\sigma^{il}_1(s,\cdot)\hat{\mu}_s, \frac{\partial (S_{\e}F)(\cdot)}{\partial x_i}>=<\sigma^{il}_1(s,\cdot)\hat{\mu}_s, S_{\e}\frac{\partial F(\cdot)}{\partial x_i}>\no\\
&=&<S_{\e}(\sigma^{il}_1(s,\cdot)\hat{\mu}_s), \frac{\partial F(\cdot)}{\partial x_i}>_{\mH}=-<\frac{\partial}{\partial x_i}S_{\e}(\sigma^{il}_1(s,\cdot)\hat{\mu}_s), F>_{\mH}.
\label{i2es}
\ee

In the following, based on Lemma \ref{prose}, we deal with $I_3, I_4$ to obtain that 
\be
&& I_3=<h^l(s,\cdot,Y_s)\hat{\mu}_s, S_{\e}F>=<S_{\e}(h^l(s,\cdot,Y_s)\hat{\mu}_s), F>_{\mH}, \label{i3es}\\
&& I_4=<(\lambda(s,\cdot,u)-1)\hat{\mu}_{s-}, S_{\e}F>=<S_{\e}((\lambda(s,\cdot,u)-1)\hat{\mu}_{s-}), F>_{\mH}.
\label{i4es}
\ee
Combining (\ref{i1es})-(\ref{i4es}) with (\ref{esex}), one can get (\ref{esexo}). 

{\bf Step 2} We prove (\ref{extrop}). 

Applying the It\^o formula to $|<Z_t^\e, F>_{\mH}|^2$, we obtain that
\ce
&&|<Z_t^\e, F>_{\mH}|^2\\
&=&|<Z_0^\e, F>_{\mH}|^2-2\int_0^t<Z_s^\e, F>_{\mH}<\frac{\partial}{\partial x_i}S_{\e}(b^i_1(s,\cdot)\hat{\mu}_s),  F>_{\mH}\dif s\no\\
&&+\int_0^t<Z_s^\e, F>_{\mH}<\frac{\partial^2}{\partial x_i\partial x_j}S_{\e}(\sigma_0^{ik}(s,\cdot) \sigma_0^{kj}(s,\cdot)\hat{\mu}_s),  F>_{\mH}\dif s\no\\
&&+\int_0^t<Z_s^\e, F>_{\mH}<\frac{\partial^2}{\partial x_i\partial x_j}
S_{\e}(\sigma_1^{il}(s,\cdot) \sigma_1^{lj}(s,\cdot)\hat{\mu}_s), F>_{\mH}\dif s\no\\
&&+2\int_0^t\int_{\mU_1}<Z_s^\e, F>_{\mH}\bigg[<S_{\e}\hat{\mu}_s, F\big(\cdot+f_1(s,\cdot,u)\big)>_{\mH}-<S_{\e}\hat{\mu}_s, F>_{\mH}\no\\
&&\qquad\qquad\qquad\qquad\qquad\qquad +<\frac{\partial}{\partial x_i}S_{\e}(f^i_1(s,\cdot,u)\hat{\mu}_s),  F>_{\mH}\bigg]\nu_1(\dif u)\dif s\no\\
&&-2\int_0^t<Z_s^\e, F>_{\mH}<\frac{\partial}{\partial x_i}S_{\e}(\sigma^{il}_1(s,\cdot)\hat{\mu}_s), F>_{\mH}\dif \hat{W}^l_s\no\\
&&+2\int_0^t<Z_s^\e, F>_{\mH}<S_{\e}(h^l(s,\cdot,Y_s)\hat{\mu}_s), F>_{\mH}\dif \hat{W}^l_s\no\\
&&+\int_0^t\int_{\mU_2}\bigg[|<Z_{s-}^\e, F>_{\mH}+<S_{\e}((\lambda(s,\cdot,u)-1)\hat{\mu}_{s-}), F>_{\mH}|^2\no\\
&&\qquad\qquad\qquad\qquad -|<Z_{s-}^\e, F>_{\mH}|^2\bigg]\tilde{\hat{N}}(\dif s, \dif u)\no\\ 
&&+\sum_{l=1}^m\int_0^t|<\frac{\partial}{\partial x_i}S_{\e}(\sigma^{il}_1(s,\cdot)\hat{\mu}_s), F>_{\mH}|^2\dif s+\sum_{l=1}^m\int_0^t|<S_{\e}(h^l(s,\cdot,Y_s)\hat{\mu}_s), F>_{\mH}|^2\dif s\no\\ 
&&-\int_0^t<\frac{\partial}{\partial x_i}S_{\e}(\sigma^{il}_1(s,\cdot)\hat{\mu}_s), F>_{\mH}<S_{\e}(h^l(s,\cdot,Y_s)\hat{\mu}_s), F>_{\mH}\dif s\no\\ 
&&+\int_0^t\int_{\mU_2}\bigg[|<Z_s^\e, F>_{\mH}+<S_{\e}((\lambda(s,\cdot,u)-1)\hat{\mu}_{s}), F>_{\mH}|^2-|<Z_s^\e, F>_{\mH}|^2\no\\ 
&&\qquad\qquad\qquad\qquad-2<Z_s^\e, F>_{\mH}<S_{\e}((\lambda(s,\cdot,u)-1)\hat{\mu}_{s}), F>_{\mH}\bigg]\nu_2(\dif u)\dif s,\no\\  
&&\qquad\qquad t\in[0,T], \quad F\in\cC_c^\infty(\mR^n).
\de
Taking $F=\phi_j$, $j=1,2,...$ and using the equality $\|Z_t^\e\|_{\mH}^2=\sum_{j=1}^{\infty}|<Z_t^\e, \phi_j>_{\mH}|^2$, we furthermore have that
\ce
\|Z_t^\e\|^2_{\mH}&=&\|Z_0^\e\|^2_{\mH}-2\int_0^t<Z_s^\e, \frac{\partial}{\partial x_i}S_{\e}(b^i_1(s,\cdot)\hat{\mu}_s)>_{\mH}\dif s\no\\
&&+\int_0^t<Z_s^\e, \frac{\partial^2}{\partial x_i\partial x_j}S_{\e}(\sigma_0^{ik}(s,\cdot) \sigma_0^{kj}(s,\cdot)\hat{\mu}_s)>_{\mH}\dif s\no\\
&&+\int_0^t<Z_s^\e,\frac{\partial^2}{\partial x_i\partial x_j}
S_{\e}(\sigma_1^{il}(s,\cdot) \sigma_1^{lj}(s,\cdot)\hat{\mu}_s)>_{\mH}\dif s\no\\
&&+2\int_0^t\int_{\mU_1}\bigg[\sum_{j=1}^{\infty}<Z_s^\e, \phi_j>_{\mH}<S_{\e}\hat{\mu}_s, \phi_j\big(\cdot+f_1(s,\cdot,u)\big)>_{\mH}-<Z_s^\e, S_{\e}\hat{\mu}_s>_{\mH}\no\\
&&\qquad\qquad\qquad\qquad +<Z_s^\e,\frac{\partial}{\partial x_i}S_{\e}(f^i_1(s,\cdot,u)\hat{\mu}_s)>_{\mH}\bigg]\nu_1(\dif u)\dif s\no\\
&&-2\int_0^t<Z_s^\e, \frac{\partial}{\partial x_i}S_{\e}(\sigma^{il}_1(s,\cdot)\hat{\mu}_s)>_{\mH}\dif \hat{W}^l_s+2\int_0^t<Z_s^\e, S_{\e}(h^l(s,\cdot,Y_s)\hat{\mu}_s)>_{\mH}\dif \hat{W}^l_s\no\\
&&+\int_0^t\int_{\mU_2}\bigg[\|Z_{s-}^\e+S_{\e}((\lambda(s,\cdot,u)-1)\hat{\mu}_{s-})\|_{\mH}^2
-\|Z_{s-}^\e\|_{\mH}^2\bigg]\tilde{\hat{N}}(\dif s, \dif u)\no\\ 
&&+\sum_{l=1}^m\int_0^t\|\frac{\partial}{\partial x_i}S_{\e}(\sigma^{il}_1(s,\cdot)\hat{\mu}_s)\|_{\mH}^2\dif s
+\sum_{l=1}^m\int_0^t\|S_{\e}(h^l(s,\cdot,Y_s)\hat{\mu}_s)\|_{\mH}^2\dif s\no\\ 
&&-\int_0^t<\frac{\partial}{\partial x_i}S_{\e}(\sigma^{il}_1(s,\cdot)\hat{\mu}_s),S_{\e}(h^l(s,\cdot,Y_s)\hat{\mu}_s)>_{\mH}\dif s\no\\ 
&&+\int_0^t\int_{\mU_2}\|S_{\e}((\lambda(s,\cdot,u)-1)\hat{\mu}_{s})\|_{\mH}^2\nu_2(\dif u)\dif s. 
\de
Thus, (\ref{extrop}) is obtained by taking the expectation under the measure $\hat{\mP}$ on two hand sides of the above equality. The proof is complete.

{\bf The proof of Lemma \ref{prozak}.}

By Lemma \ref{trex}, it holds that 
\be
\hat{\mE}\|Z^\e_t\|^2_{\mH}&=&\|Z^\e_0\|^2_{\mH}+I_1+I_2+I_3+I_4+I_5+I_6+I_7+I_8,
\label{extrop1}
\ee
where 
\ce
&&I_1:=-2\int_0^t\hat{\mE}<Z^{\e}_s, \frac{\partial}{\partial x_i}S_{\e}(b^i_1(s,\cdot)\hat{\mu}_s)>_{\mH}\dif s,\\
&&I_2:=\int_0^t\hat{\mE}<Z^{\e}_s, \frac{\partial^2}{\partial x_i\partial x_j}S_{\e}(\sigma_0^{ik}(s,\cdot) \sigma_0^{kj}(s,\cdot)\hat{\mu}_s)>_{\mH}\dif s,\\
&&I_3:=\int_0^t\hat{\mE}<Z^{\e}_s,\frac{\partial^2}{\partial x_i\partial x_j}
S_{\e}(\sigma_1^{il}(s,\cdot) \sigma_1^{lj}(s,\cdot)\hat{\mu}_s)>_{\mH}\dif s,\\
&&I_4:=2\int_0^t\int_{\mU_1}\hat{\mE}\bigg[\sum_{j=1}^{\infty}<Z^{\e}_s, \phi_j>_{\mH}<Z^{\e}_s, \phi_j\big(\cdot+f_1(s,\cdot,u)\big)>_{\mH}-\|Z^{\e}_s\|^2_{\mH}\no\\
&&\qquad\qquad\qquad\qquad\qquad+<Z^{\e}_s,\frac{\partial}{\partial x_i}S_{\e}(f^i_1(s,\cdot,u)\hat{\mu}_s)>_{\mH}\bigg]\nu_1(\dif u)\dif s,\no\\
&&I_5:=\sum_{l=1}^m\int_0^t\hat{\mE}\|\frac{\partial}{\partial x_i}S_{\e}(\sigma^{il}_1(s,\cdot)\hat{\mu}_s)\|_{\mH}^2\dif s,\\
&&I_6:=\sum_{l=1}^m\int_0^t\hat{\mE}\|S_{\e}(h^l(s,\cdot,Y_s)\hat{\mu}_s)\|_{\mH}^2\dif s,\\
&&I_7:=-\int_0^t\hat{\mE}<\frac{\partial}{\partial x_i}S_{\e}(\sigma^{il}_1(s,\cdot)\hat{\mu}_s),S_{\e}(h^l(s,\cdot,Y_s)\hat{\mu}_s)>_{\mH}\dif s,\\
&&I_8:=\int_0^t\int_{\mU_2}\hat{\mE}\|S_{\e}((\lambda(s,\cdot,u)-1)\hat{\mu}_s)\|_{\mH}^2\nu_2(\dif u)\dif s.
\de
By Lemma \ref{onede}, we know that 
\be\label{17es}
I_1+I_6+I_7\leq C\int_0^t\hat{\mE}\|Z^{\e}_s\|^2_{\mH}\dif s.
\ee
And then Lemma \ref{secde} admits us to obtain that 
\be\label{235es}
I_2+I_3+I_5\leq C\int_0^t\hat{\mE}\|Z^{\e}_s\|^2_{\mH}\dif s.
\ee

To estimate $I_4$, we divide $I_4$ into three parts $I_{41},  I_{42}, I_{43}$. For $I_{41}$, since for $s\in[0,T]$ and $u\in\mU_1$, $x\to y:=x+f_1(s,x,u)$ is invertible, there exists an inverse function $x=g(s,y,u)$. Moreover, 
$$
\left|\det(\partial_y g(s,y,u))\right|=\left|\det\((J_{f_1}+I)^{-1}\)\right|=|\det(J_{f_1}+I)|^{-1}\leq G_3(u).
$$
And then for any $j=1, 2, \cdots$, it holds that 
\ce
<Z^{\e}_s, \phi_j\big(\cdot+f_1(s,\cdot,u)\big)>_{\mH}&=&\int_{\mR^n}Z^{\e}_s(x)\phi_j\big(x+f_1(s,x,u)\big)\dif x\\
&=&\int_{\mR^n}Z^{\e}_s\(g(s,y,u)\)\phi_j(y)|\det(\partial_y g(s,y,u))|\dif y\\
&=&<Z^{\e}_s\(g(s,\cdot,u)\)|\det(\partial_y g(s,\cdot,u))|, \phi_j>_{\mH}.
\de
Thus, by the above deduction, we have that
\be
I_{41}&=&2\int_0^t\int_{\mU_1}\hat{\mE}\sum_{j=1}^{\infty}<Z^{\e}_s, \phi_j>_{\mH}<Z^{\e}_s\(g(s,\cdot,u)\)|\det(\partial_y g(s,\cdot,u))|, \phi_j>_{\mH}\nu_1(\dif u)\dif s\no\\
&=&2\int_0^t\int_{\mU_1}\hat{\mE}<Z^{\e}_s\(g(s,\cdot,u)\)|\det(\partial_y g(s,\cdot,u))|, Z^{\e}_s>_{\mH}\nu_1(\dif u)\dif s\no\\
&\leq&2\int_0^t\int_{\mU_1}\hat{\mE}\|Z^{\e}_s\|_{\mH}\left\|Z^{\e}_s\(g(s,\cdot,u)\)|\det(\partial_y g(s,\cdot,u))|\right\|_{\mH}\nu_1(\dif u)\dif s\no\\
&\leq&\int_0^t\int_{\mU_1}\left(\hat{\mE}\|Z^{\e}_s\|^2_{\mH}+\hat{\mE}\left\|Z^{\e}_s\(g(s,\cdot,u)\)\det(\partial_y g(s,\cdot,u))\right\|^2_{\mH}\right)\nu_1(\dif u)\dif s.
\label{i411}
\ee
Note that 
\be
\left\|Z^{\e}_s\(g(s,\cdot,u)\)\det(\partial_y g(s,\cdot,u))\right\|^2_{\mH}&=&\int_{\mR^n}\left|Z^{\e}_s\(g(s,y,u)\)\right|^2|\det(\partial_y g(s,y,u))|^2\dif y\no\\
&=&\int_{\mR^n}\left|Z^{\e}_s\(x\)\right|^2|\det(\partial_y g(s,y,u))|\dif x\no\\
&\leq&G_3(u)\int_{\mR^n}\left|Z^{\e}_s\(x\)\right|^2\dif x=G_3(u)\|Z^{\e}_s\|^2_{\mH}.
\label{varichan}
\ee
So, inserting (\ref{varichan}) into (\ref{i411}), we furthermore obtain that
\be
I_{41}\leq \int_0^t\int_{\mU_1}(G_3(u)+1)\hat{\mE}\|Z^{\e}_s\|^2_{\mH}\nu_1(\dif u)\dif s.
\label{i412}
\ee 

For $I_{43}$, note that 
\ce
&&2<Z^{\e}_s,\frac{\partial}{\partial x_i}S_{\e}(f^i_1(s,\cdot,u)\hat{\mu}_s)>_{\mH}\\
&=&2<Z^{\e}_s,f^i_1(s,\cdot,u)\frac{\partial}{\partial x_i}Z^{\e}_s>_{\mH}+2<Z^{\e}_s,\frac{\partial}{\partial x_i}S_{\e}(f^i_1(s,\cdot,u)\hat{\mu}_s)-f^i_1(s,\cdot,u)\frac{\partial}{\partial x_i}Z^{\e}_s>_{\mH}\\
&\leq&2\int_{\mR^n}Z^{\e}_s(x)f^i_1(s,x,u)\frac{\partial}{\partial x_i}Z^{\e}_s(x)\dif x+\|Z^{\e}_s\|^2_{\mH}+\left\|\frac{\partial}{\partial x_i}S_{\e}(f^i_1(s,\cdot,u)\hat{\mu}_s)-f^i_1(s,\cdot,u)\frac{\partial}{\partial x_i}Z^{\e}_s\right\|^2_{\mH}\\
&=&\int_{\mR^n}f^i_1(s,x,u)\frac{\partial}{\partial x_i}|Z^{\e}_s(x)|^2\dif x+\|Z^{\e}_s\|^2_{\mH}+\left\|\frac{\partial}{\partial x_i}S_{\e}(f^i_1(s,\cdot,u)\hat{\mu}_s)-f^i_1(s,\cdot,u)\frac{\partial}{\partial x_i}Z^{\e}_s\right\|^2_{\mH}\\
&\leq&L^\prime_1(s)G_1(u)\|Z^{\e}_s\|^2_{\mH}+\|Z^{\e}_s\|^2_{\mH}+\left\|\frac{\partial}{\partial x_i}S_{\e}(f^i_1(s,\cdot,u)\hat{\mu}_s)-f^i_1(s,\cdot,u)\frac{\partial}{\partial x_i}Z^{\e}_s\right\|^2_{\mH},
\de
where we use the formula for integration by parts, ($\mathbf{H}^3_{f_1}$) ($\mathbf{H}^{1'}_{b_1, \sigma_0, \sigma_1, f_1}$) and the fact that $\lim\limits_{|x|\to\infty}Z^{\e}_s(x)=0$ in the last equality. So, for the third term in the right side of the above inequality, it follows from the definition of $S_{\e}$ that
\ce
&&\left\|\frac{\partial}{\partial x_i}S_{\e}(f^i_1(s,\cdot,u)\hat{\mu}_s)-f^i_1(s,\cdot,u)\frac{\partial}{\partial x_i}Z^{\e}_s\right\|^2_{\mH}\\
&=&\int_{\mR^n}\left|\frac{\partial}{\partial x_i}S_{\e}(f^i_1(s,\cdot,u)\hat{\mu}_s)(x)-f^i_1(s,x,u)\frac{\partial}{\partial x_i}Z^{\e}_s(x)\right|^2\dif x\\
&=&\int_{\mR^n}\left|\int_{\mR^n}(2\pi \e)^{-\frac{n}{2}}(f^i_1(s,y,u)-f^i_1(s,x,u))\frac{\partial}{\partial x_i}\left(\exp\left\{-\frac{|x-y|^2}{2\e}\right\}\right)\hat{\mu}_s(\dif y)\right|^2\dif x\\
&\leq&\int_{\mR^n}\left(\int_{\mR^n}(2\pi \e)^{-\frac{n}{2}}|f^i_1(s,y,u)-f^i_1(s,x,u)|\frac{|x-y|}{\e}\exp\left\{-\frac{|x-y|^2}{2\e}\right\}\hat{\mu}_s(\dif y)\right)^2\dif x\\
&\leq&L^\prime_1(s)G^2_1(u)\int_{\mR^n}\left(\int_{\mR^n}(2\pi \e)^{-\frac{n}{2}}\frac{|x-y|^2}{\e}\exp\left\{-\frac{|x-y|^2}{2\e}\right\}\hat{\mu}_s(\dif y)\right)^2\dif x\\
&=&L^\prime_1(s)G^2_1(u)\int_{\mR^n}\bigg(\int_{\mR^n}2^{\frac{n}{2}}\frac{|x-y|^2}{\e}\exp\left\{-\frac{|x-y|^2}{4\e}\right\}\\
&&\qquad\qquad \times(2\pi 2\e)^{-\frac{n}{2}}\exp\left\{-\frac{|x-y|^2}{4\e}\right\}\hat{\mu}_s(\dif y)\bigg)^2\dif x\\
&\overset{ve^{-\frac{v}{4}}\leq C, v\geq 0}{\leq}&CL^\prime_1(s)G^2_1(u)\int_{\mR^n}|S_{2\e}\hat{\mu}_s(x)|^2\dif x=CL^\prime_1(s)G^2_1(u)\|S_{2\e}\hat{\mu}_s\|^2_{\mH}\\
&\leq& CL^\prime_1(s)G^2_1(u)\|Z^{\e}_s\|^2_{\mH},
\de
where the last inequality is based on Lemma \ref{prose}. Thus, by ($\mathbf{H}^{1'}_{b_1, \sigma_0, \sigma_1, f_1}$) we have that
\be\label{4es}
I_4&\leq&\int_0^t\int_{\mU_1}\left(G_3(u)\hat{\mE}\|Z^{\e}_s\|^2_{\mH}+L^\prime_1(s)G_1(u)\hat{\mE}\|Z^{\e}_s\|^2_{\mH}+CL^\prime_1(s)G^2_1(u)\hat{\mE}\|Z^{\e}_s\|^2_{\mH}\right)\nu_1(\dif u)\dif s\no\\
&\leq&C\int_0^t\hat{\mE}\|Z^{\e}_s\|^2_{\mH}\dif s. 
\ee

For $I_8$, we rewrite it to obtain that 
\ce
&&\int_{\mU_2}\|S_{\e}((\lambda(s,\cdot,u)-1)\hat{\mu}_{s})\|_{\mH}^2\nu_2(\dif u)
=\int_{\mU_2}\nu_2(\dif u)\int_{\mR^n}|S_{\e}((\lambda(s,\cdot,u)-1)\hat{\mu}_{s})(x)|^2\dif x\\
&\leq&\int_{\mU_2}\nu_2(\dif u)\int_{\mR^n}\dif x\int_{\mR^n}(2\pi \e)^{-\frac{n}{2}}\exp\left\{-\frac{|x-y|^2}{2\e}\right\}(\lambda(s,y,u)-1)\hat{\mu}_{s}(\dif y)\\
&&\qquad\qquad \times\int_{\mR^n}(2\pi \e)^{-\frac{n}{2}}\exp\left\{-\frac{|x-z|^2}{2\e}\right\}(\lambda(s,z,u)-1)\hat{\mu}_{s}(\dif z)\\
&=&\int_{\mR^n}\dif x\int_{\mR^n}(2\pi \e)^{-\frac{n}{2}}\exp\left\{-\frac{|x-y|^2}{2\e}\right\}\hat{\mu}_{s}(\dif y)\int_{\mR^n}(2\pi \e)^{-\frac{n}{2}}\exp\left\{-\frac{|x-z|^2}{2\e}\right\}\hat{\mu}_{s}(\dif z)\\
&&\qquad\qquad \times\int_{\mU_2}(\lambda(s,y,u)-1)(\lambda(s,z,u)-1)\nu_2(\dif u)\\
&\overset{(\mathbf{H}_{\lambda})}{\leq}&C\int_{\mR^n}\dif x\int_{\mR^n}(2\pi \e)^{-\frac{n}{2}}\exp\left\{-\frac{|x-y|^2}{2\e}\right\}\hat{\mu}_{s}(\dif y)\int_{\mR^n}(2\pi \e)^{-\frac{n}{2}}\exp\left\{-\frac{|x-z|^2}{2\e}\right\}\hat{\mu}_{s}(\dif z)\\
&=&C\int_{\mR^n}|S_{\e}\hat{\mu}_{s}(x)|^2\dif x=C\|Z^{\e}_s\|^2_{\mH},
\de
and then
\be\label{8es}
I_8\leq C\int_0^t\hat{\mE}\|Z^{\e}_s\|^2_{\mH}\dif s.
\ee

By combining (\ref{17es}) (\ref{235es}) (\ref{4es}) (\ref{8es}) with (\ref{extrop1}), it holds that
\ce
\hat{\mE}\|Z^{\e}_t\|^2_{\mH}\leq\|Z^\e_0\|^2_{\mH}+C\int_0^t\hat{\mE}\|Z^{\e}_s\|^2_{\mH}\dif s.
\de
This is just right (\ref{moin}).

In the following, we prove $\hat{\mu}_t\in\mH, a.s.$ for $t\in[0,T]$. By the above inequality and the Gronwall inequality, we have that
\ce
\hat{\mE}\|Z^{\e}_t\|^2_{\mH}\leq\|Z^\e_0\|^2_{\mH}e^{Ct}.
\de
Thus, it follows from the Fatou lemma that
\ce
&&\hat{\mE}\|\hat{\mu}_t\|^2_{\mH}=\hat{\mE}\left(\sum_{j=1}^{\infty}<\phi_j, \hat{\mu}_t>^2\right)=\hat{\mE}\left(\sum_{j=1}^{\infty}\lim_{\e\rightarrow 0}(<S_\e\phi_j, \hat{\mu}_t>)^2\right)\\
&\leq&\liminf_{\e\rightarrow 0}\hat{\mE}\left(\sum_{j=1}^{\infty}<S_\e\phi_j, \hat{\mu}_t>^2\right)=\liminf_{\e\rightarrow 0}\hat{\mE}\left(\sum_{j=1}^{\infty}<\phi_j, S_\e\hat{\mu}_t>_{\mH}^2\right)\\
&=&\liminf_{\e\rightarrow 0}\hat{\mE}\|Z^{\e}_t\|^2_{\mH}\leq \liminf_{\e\rightarrow 0}\|Z^\e_0\|^2_{\mH}e^{Ct}=\|\hat{\mu}_0\|^2_{\mH}e^{Ct}<\infty.
\de
That is, $\hat{\mu}_t\in\mH, a.s.$ for $t\in[0,T]$.

 \bigskip

\textbf{Acknowledgements:}

The author would like to thank Professor Xicheng Zhang for valuable discussions and also thanks Professor Renming Song for providing her an excellent environment to work in the University of Illinois at Urbana-Champaign. Besides, the author is 
grateful to two referees since their suggestions and comments allowed her to improve the results and the presentation of this paper.

\end{document}